\documentclass[9pt,11pt]{article}%
\usepackage{amsfonts}
\usepackage[applemac]{inputenc}
\usepackage{amssymb}
\usepackage{makeidx}
\usepackage{graphicx}
\usepackage{amsmath}
\usepackage[a4paper]{geometry}
\usepackage[displaymath]{lineno}
\usepackage[singlespacing]{setspace}
\usepackage[normalem]{ulem}
\usepackage{lineno}
\usepackage{color}%
\providecommand{\U}[1]{\protect\rule{.1in}{.1in}}
\providecommand{\U}[1]{\protect\rule{.1in}{.1in}}

\newtheorem{theo}{Theorem}
\newtheorem{lem}[theo]{Lemma}

\newtheorem{cor}[theo]{Corollary}
\newtheorem{rem}{Remark}

\newenvironment{dem}[1][Proof]{\noindent \textbf{#1.} }{\ \rule{0.5em}{0.5em}}
\begin{document}

\title{Valadier-like formulas for the supremum function II: The compactly indexed
case\thanks{Research of the first and the second authors is supported by
CONICYT grants, Fondecyt no. 1150909 and 1151003, and Basal PFB-03. Research
of the second and third authors is supported by MINECO of Spain and FEDER of
EU, grant MTM2014-59179-C2-1-P. Research of the third author is also supported
by the Australian Research Council, Project DP160100854.}}
\author{R. Correa\thanks{e-mail: rcorrea@dim.uchile.cl}, A.\ Hantoute\thanks{e-mail:
ahantoute@dim.uchile.cl (corresponding author)}, M.A. López\thanks{e-mail:
marco.antonio@ua.es}\\${^{\dag}}${\scriptsize Universidad de O'Higgins, Chile, and DIM-CMM of
Universidad de Chile}\\$^{\ddag}${\scriptsize Center for Mathematical Modeling (CMM), Universidad de
Chile }\\$^{§}${\scriptsize Universidad de Alicante, Spain, and CIAO, Federation
University, Ballarat, Australia}}
\date{}
\maketitle

\begin{abstract}
Continuing with the work on the subdifferential of the pointwise supremum of
convex functions, started in \emph{Valadier-like formulas for the supremum
function I }\cite{RcAhMl16}, we focus now on the compactly indexed case. We
assume that the index set is compact and that the data functions are upper
semicontinuous with respect to the index variable (actually, this assumption
will only affect the set of $\varepsilon$-active indices at the reference
point). As in the previous work, we do not require any continuity assumption
with respect to the decision variable. The current compact setting gives rise
to more explicit formulas, which only involve subdifferentials at the
reference point of active data functions. Other formulas are derived under
weak continuity assumptions. These formulas reduce to the characterization
given by Valadier \cite[Theorem 2]{Valadier} when the supremum function is continuous.



\end{abstract}

\section{Introduction\label{Sect1}}

Let us consider\ a family of convex functions$\ f_{t}:X\rightarrow
\overline{\mathbb{R}}:=\mathbb{R}\cup\{\pm\infty\},$ $t\in T,$ defined in a
locally convex topological vector space $X.$ The aim of this paper, which
continues \cite{RcAhMl16},\ is to\ give characterizations of the
subdifferential of the supremum function
\begin{equation}
f:=\sup_{t\in T}f_{t}, \label{suprefunction}%
\end{equation}
which only involve the exact subdifferentials of data functions at the
reference point rather than at nearby ones. Our results are based on some
compactness assumptions of certain subsets of the index set $T$, and some
upper semicontinuity assumptions of the mappings $t\mapsto f_{t}(z),$
$z\in\operatorname*{dom}f$. In Theorem \ref{thmcompact} we establish that,
under a natural closure condition (\ref{assumption}), for all $x\in X$%
\[
\partial f(x)=\bigcap_{L\in\mathcal{F}(x)}\overline{\operatorname*{co}%
}\left\{  \bigcup\limits_{t\in T(x)}\partial(f_{t}+\mathrm{I}_{L\cap
\operatorname*{dom}f})(x)\right\}  ,
\]
where
\[
T(x):=\{t\in T\mid f_{t}(x)=f(x)\},
\]
and
\[
\mathcal{F}(x):=\{\text{finite-dimensional linear subspaces of }X\text{ such
that }x\in L\}.
\]
The following characterization of $\partial f(x),$ given\ in Theorem
\ref{spe1}, uses the $\varepsilon$-subdifferentials of the functions $f_{t}:$%
\[
\partial f(x)=\bigcap_{\substack{\varepsilon>0\\L\in\mathcal{F}(x)}%
}\overline{\operatorname*{co}}\left\{  \bigcup\limits_{t\in T(x)}%
\partial_{\varepsilon}f_{t}(x)+\mathrm{N}_{L\cap\operatorname*{dom}%
f}(x)\right\}  .
\]

Condition (\ref{assumption}) covers the case when the functions $f_{t}$ are
lower semicontinuous (lsc). In particular, if the restriction of $f$ to the
affine hull of $\operatorname*{dom}f$ is continuous on the relative interior
of $\operatorname*{dom}f$ (assumed to be nonempty), the intersection over
$L\in\mathcal{F}(x)$ can be removed, giving rise to the following formula,
established in Theorem \ref{thmcompact0}:%
\[
\partial f(x)=\overline{\operatorname*{co}}\left\{  \bigcup\limits_{t\in
T(x)}\partial(f_{t}+\mathrm{I}_{\operatorname*{dom}f})(x)\right\}  .
\]

Our results generalize and improve the well-known\textbf{\ }formula\ due to
Valadier\ \cite[Theorem 2]{Valadier}, which establishes that, under the
continuity of $f$ at $x,$
\begin{equation}
\partial f(x)=\overline{\operatorname*{co}}\left\{  \bigcup\limits_{t\in
T(x)}\partial f_{t}(x)\right\}  . \label{Valadier-like formula}%
\end{equation}
Actually, we show that if $f$ is continuous at some point (not necessarily the
reference point $x),\,$then we get
\[
\partial f(x)=\mathrm{N}_{\operatorname*{dom}f}(x)+\overline
{\operatorname*{co}}\left\{  \bigcup\limits_{t\in T(x)}\partial f_{t}%
(x)\right\}  ,
\]
which reduces to (\ref{Valadier-like formula})$\ $whenever $\mathrm{N}%
_{\operatorname*{dom}f}(x)=\left\{  \theta\right\}  .$ Observe that the
continuity $f$ at the reference point $x$ is equivalent to the continuity of
$f$ at some point together with $\mathrm{N}_{\operatorname*{dom}f}(x)=\left\{
\theta\right\}  .$ Indeed, the last condition implies that $x$ is in the
\emph{quasi-interior} of $\operatorname*{dom}f$ which coincides with its
interior (see, i.e., \cite{Borwlewis}).

There is a wide literature dealing with subdifferential calculus rules for the
supremum of convex functions; we refer for instance to \cite{Brondsted72,
Han06, HanLop08, HanLopZal2008, Ioffe11, Ioffe72, LiNg11, Soloviev01,
Valadier, ZalinescuBook}, among many others. We also refer to \cite{MOR}, and
references therein, for the nonconvex case. The supremum function plays a
crucial role in many fields, including semi-infinite optimization
(\cite{WELL}, \cite{Zheng}). 

The paper is organized as follows. After Section 2, devoted to preliminaries,
main Section 3 provides the desired characterization of $\partial f(x)$ in
different settings: Theorem \ref{thmcompact0} deals with a
finite-dimensional-like setting, where function $f$ satisfies a weak
continuity condition, which is always held in finite dimensions; Theorem
\ref{thmcompact1} concerns the supremum of lsc convex functions, while the
most general result is given in Theorem \ref{thmcompact} under a closure-type
condition. All these results use the exact subdifferential of the data
functions at the nominal point. Another formula using approximate
subdifferentials of data functions is given in Theorem \ref{spe1}. Finally,
Theorem \ref{corcompcont} provides a simpler formula similar to
(\ref{Valadier-like formula})$\ $when additional continuity assumptions are imposed.

\section{Preliminaries\label{Sect2}}

In this paper $X$ stands for a (real) separated locally convex (lcs, shortly)
space\textbf{,} whose\textbf{\ }topological dual space is denoted by $X^{\ast
}$ and endowed with the weak*-topology.$\ $Hence, $X$ and $X^{\ast}$ form a
dual pair by means of the canonical bilinear form $\langle x,x^{\ast}%
\rangle=\langle x^{\ast},x\rangle:=x^{\ast}(x),$ $(x,x^{\ast})\in X\times
X^{\ast}.$ The zero vectors are denoted by $\theta$, and the convex, closed
and balanced neighborhoods of $\theta$ are called $\theta$-neighborhoods. The
family of such $\theta$-neighborhoods in $X$ and in $X^{\ast}$ are denoted by
$\mathcal{N}_{X}$ and $\mathcal{N}_{X^{\ast}},$ respectively.

Given a nonempty set $A$ in $X$ (or in $X^{\ast}$), by $\operatorname*{co}A$
and $\operatorname*{aff}A$ we denote the \emph{convex hull} and\emph{\ }the
\emph{affine hull }of $A$, respectively. Moreover, $\operatorname*{cl}A$ and
$\overline{A}$ are indistinctly used for denoting the \emph{closure }of $A$
(\emph{weak*-}$\emph{closure}$ if $A\subset X^{\ast}$). Thus, $\overline
{\operatorname*{co}}A:=\operatorname*{cl}(\operatorname*{co}A)$,
$\overline{\operatorname*{aff}}A:=\operatorname*{cl}(\operatorname*{aff}A),$
etc. We use $\operatorname*{ri}A$ to denote the (topological) \emph{relative
interior }of $A$ (i.e., the interior of $A$ in the topology relative to
$\operatorname*{aff}A$ when this set is closed, and the empty set otherwise).
We consider the \emph{orthogonal }of\emph{\ }$A$ defined by
\[
A^{\perp}:=\left\{  x^{\ast}\in X^{\ast}\mid\langle x^{\ast},x\rangle=0\text{
for all }x\in A\right\}  .
\]

We say that a convex function $\varphi:X\rightarrow\overline{\mathbb{R}}$ is
proper if its \emph{(effective)} \emph{domain,} $\operatorname*{dom}%
\varphi:=\{x\in X\mid\varphi(x)<+\infty\}$, is nonempty and it does not take
the value $-\infty.$ The \emph{lsc envelope }of $\varphi$ is denoted by
$\operatorname*{cl}\varphi$.\emph{\ }We adopt the convention\emph{ }$\left(
+\infty\right)  +(-\infty)=\left(  -\infty\right)  +(+\infty)=+\infty.$

If $\psi:X\rightarrow\mathbb{R\cup}\left\{  +\infty\right\}  $ is another
proper convex function, which is finite and continuous at some point in
$\operatorname*{dom}\varphi,$ then we have \cite[Corollary 9(iii)]{HanLop08}%
\begin{equation}
\operatorname*{cl}(\max\left\{  \varphi,\psi\right\}  )=\max\left\{
\operatorname*{cl}\varphi,\operatorname*{cl}\psi\right\}  . \label{cl}%
\end{equation}

For $\varepsilon\geq0$, the $\varepsilon$-\emph{subdifferential} of $\varphi$
at a point $x$ where $\varphi(x)$ is finite is the weak*-closed convex set
\[
\partial_{\varepsilon}\varphi(x):=\{x^{\ast}\in X^{\ast}\mid\varphi
(y)-\varphi(x)\geq\langle x^{\ast},y-x\rangle-\varepsilon\text{ for all }y\in
X\}.
\]
If $\varphi(x)\notin\mathbb{R}$, then we set $\partial_{\varepsilon}%
\varphi(x)=\emptyset$. In particular, for $\varepsilon=0$ we get the
\emph{Fenchel subdifferential} of $\varphi$ at $x,$ $\partial\varphi
(x):=\partial_{0}\varphi(x)$. When $x\in\operatorname*{dom}(\partial
\varphi):=\left\{  y\in X\mid\partial\varphi(y)\neq\emptyset\right\}  ,$ we
know that
\begin{equation}
\varphi(x)=(\operatorname*{cl}\varphi)(x)\text{ and }\partial_{\varepsilon
}\varphi(x)=\partial_{\varepsilon}(\operatorname*{cl}\varphi)(x),\text{ for
all }\varepsilon\geq0\text{.}\label{noemptys}%
\end{equation}

The \emph{indicator }and the \emph{support }functions of $A\subset X$ are,
respectively, defined as
\[
\mathrm{I}_{A}(x):=0,\text{ if }x\in A;\text{ }+\infty,\text{ if }x\in
X\setminus A,
\]%
\begin{equation}
\sigma_{A}(x^{\ast}):=\sup\{\langle x^{\ast},a\rangle\mid a\in A\},\text{
}x^{\ast}\in X^{\ast}, \label{support}%
\end{equation}
with the convention $\sigma_{\emptyset}\equiv-\infty.$

If $A$ is convex and\ $x\in X$, we define the \emph{normal cone }to $A$\ at
$x$ as
\[
\mathrm{N}_{A}(x):=\{x^{\ast}\in X^{\ast}\mid\langle x^{\ast},y-x\rangle
\leq0\text{ for all }y\in A\},\text{ if \ }x\in A,
\]
and $\mathrm{N}_{A}(x)=\emptyset,$ if $x\in X\setminus A$.

Now, we review the results given in \cite{RcAhMl16}, which constitute the main
foundations of the present work. We have proved there that the subdifferential
of the supremum function $f=\sup_{t\in T}f_{t}$ is expressed in terms of
appropriate enlargements of the Fenchel subdifferential, $\breve{\partial}%
_{p}^{\varepsilon}(f_{t}+\mathrm{I}_{\overline{L\cap\operatorname*{dom}f}}),$
$t\in T_{\varepsilon}(x),$ $p\in\mathcal{P}$, and $L\in\mathcal{F}(x),$ where
\[
T_{\varepsilon}(x):=\left\{  t\in T\mid f_{t}(x)\geq f(x)-\varepsilon\right\}
,\text{ }\varepsilon\geq0,
\]%
\[
\mathcal{P}=\{\text{continuous seminorms on }X\},
\]
and
\[
\mathcal{F}(x):=\{\text{finite-dimensional linear subspaces }L\subset X\text{
containing }x\}.
\]
Such enlargements involve the exact subdifferentials of functions
$f_{t}+\mathrm{I}_{\overline{L\cap\operatorname*{dom}f}}$ at nearby points of
$x.$\ Precisely, for a convex function $\varphi$, $\breve{\partial}%
_{p}^{\varepsilon}\varphi(x)$ is defined\ by%
\[
\breve{\partial}_{p}^{\varepsilon}\varphi(x):=\left\{  y^{\ast}\in
\partial\varphi(y)\mid p(y-x)\leq\varepsilon,\text{ }\left\vert \varphi
(y)-\varphi(x)\right\vert \leq\varepsilon,\text{ and }\left\vert \langle
y^{\ast},y-x\rangle\right\vert \leq\varepsilon\right\}  .
\]
Observe that $\breve{\partial}_{p}^{\varepsilon}\varphi$ provides\ an outer
approximation of $\partial\varphi$ as
\begin{equation}
\partial\varphi(x)\subset\breve{\partial}_{p}^{\varepsilon}\varphi
(x)\subset\partial_{2\varepsilon}\varphi(x).\label{rel}%
\end{equation}

When the functions $f_{t}$ are proper and lsc, we proved in \cite[Theorem
8]{RcAhMl16} that
\begin{equation}
\partial f(x)=\bigcap_{_{\substack{\varepsilon>0,\text{ }p\in\mathcal{P}%
\\L\in\mathcal{F}(x)}}}\overline{\operatorname*{co}}\left\{  \bigcup
\limits_{t\in T_{\varepsilon}(x)}\breve{\partial}_{p}^{\varepsilon}%
(f_{t}+\mathrm{I}_{\overline{L\cap\operatorname*{dom}f}})(x)\right\}  .
\label{mn}%
\end{equation}
In particular, if the restriction of $f$ to $\operatorname*{aff}%
(\operatorname*{dom}f)$ is continuous on $\operatorname*{ri}%
(\operatorname*{dom}f)$ (assumed to be nonempty), the intersection over
$L\in\mathcal{F}(x)$ can be removed to obtain \cite[Theorem 9]{RcAhMl16}:
\begin{equation}
\partial f(x)=\bigcap_{_{\varepsilon>0,\text{ }p\in\mathcal{P}}}%
\overline{\operatorname*{co}}\left\{  \bigcup\limits_{t\in T_{\varepsilon}%
(x)}\breve{\partial}_{p}^{\varepsilon}(f_{t}+\mathrm{I}_{\overline
{\operatorname*{dom}f}})(x)\right\}  . \label{fri}%
\end{equation}
Moreover, if the $f_{t}$'s are proper but not necessarily lsc, and $f$ is
finite and continuous at some point, we obtained\ from (\ref{fri}) that
\cite[Theorem 10]{RcAhMl16}:
\[
\partial f(x)=\mathrm{N}_{\operatorname*{dom}f}(x)+\bigcap_{\varepsilon
>0,\text{ }p\in\mathcal{P}}\overline{\operatorname*{co}}\left\{
\bigcup\limits_{t\in\tilde{T}_{\varepsilon}(x)}\breve{\partial}_{p}%
^{\varepsilon}(\operatorname*{cl}f_{t})(x)\right\}  ,
\]
where $\tilde{T}_{\varepsilon}(x):=\{t\in T\mid(\operatorname*{cl}%
f_{t})(x)\geq f(x)-\varepsilon\}.$

In the particular case when $f$ is continuous at the point $x,$ we recover in
\cite[Corollary 12]{RcAhMl16} the Valadier formula \cite[Theorem 1]%
{Valadier}:
\[
\partial f(x)=\bigcap_{\varepsilon>0,\text{ }p\in\mathcal{P}}\overline
{\operatorname*{co}}\left\{  \bigcup\limits_{t\in T_{\varepsilon}(x),\text{
}p(y-x)\leq\varepsilon}\partial f_{t}(y)\right\}  .
\]

\section{Compactly indexed case}

In this section we characterize the subdifferential of the supremum function\
\[
f=\sup_{t\in T}f_{t}\
\]
of a compactly indexed family\ of convex functions $f_{t}:X\rightarrow
\mathbb{R}\cup\left\{  \pm\infty\right\}  ,$ $t\in T,$ where $X$ is a lcs
space whose family of continuous seminorms is denoted by $\mathcal{P},$ and
the dual space $X^{\ast}$ is endowed with the weak*-topology.

First, we state our result in a finite dimensional-like setting. Recall that
\[
T_{\varepsilon}(x)=\left\{  t\in T\mid f_{t}(x)\geq f(x)-\varepsilon\right\}
,\text{ }\varepsilon\geq0,
\]
and
\[
T(x)=\left\{  t\in T\mid f_{t}(x)=f(x)\right\}  .
\]

\begin{theo}
\label{thmcompact0}Assume that the family of convex functions $\left\{
f_{t},t\in T\right\}  $ is such that the function $f_{\mid\operatorname*{aff}%
(\operatorname*{dom}f)}$ is finite and continuous on $\operatorname*{ri}%
(\operatorname*{dom}f),$ assumed to be nonempty. Suppose that
\begin{equation}
\operatorname*{cl}f=\sup_{t\in T}(\operatorname*{cl}f_{t}).\label{assumption}%
\end{equation}
Let $x\in X$ be such that for some $\varepsilon_{0}>0$\emph{:}

$(i)$ the set $T_{\varepsilon_{0}}(x)$ is compact,

$(ii)$ the functions $t\mapsto f_{t}(z),$ $z\in\operatorname*{dom}f,$\ are
upper semicontinuous (usc) on $T_{\varepsilon_{0}}(x).$ \newline Then
\begin{equation}
\partial f(x)=\overline{\operatorname*{co}}\left\{  \bigcup\limits_{t\in
T(x)}\partial(f_{t}+\mathrm{I}_{\operatorname*{dom}f})(x)\right\}  .
\label{mainf}%
\end{equation}

\end{theo}

The proof of this theorem uses the following technical lemma:

\begin{lem}
\label{lt}Given the family of convex functions $\left\{  f_{t},t\in T\right\}
$ and $x_{0}\in\operatorname*{dom}(\partial f),$ consider the functions
\[
\ell_{t}:=\max\{f_{t},f(x_{0})-c\},\text{ }t\in T,
\]
where $c>0,$ and
\[
\ell:=\sup_{t\in T}\ell_{t}.
\]
Then, under condition \emph{(\ref{assumption})}, there exists an open
neighborhood $U$ of $x_{0}$ such that the proper functions $\ell_{t}$ satisfy
the following\emph{:}

\emph{(i)} $\operatorname*{cl}f=\operatorname*{cl}\ell$ and $f=\max\left\{
f,f(x_{0})-c\right\}  =\ell,$ on $U.$

\emph{(ii) }$\operatorname*{cl}\ell=\sup_{t\in T}(\operatorname*{cl}\ell
_{t}).$

\emph{(iii) }$\{t\in T\mid\ell_{t}(x_{0})\geq\ell(x_{0})-\varepsilon
\}=T_{\varepsilon}(x_{0})$ \ for all $\varepsilon\in\left[  0,c\right[  $ $.$

\emph{(iv)} $\partial f(x_{0})=\partial\ell(x_{0}).$
\end{lem}

\begin{dem}
We may suppose that $x_{0}=\theta$ and $f(\theta)=0.$ Observe that $\ell
=\max\left\{  f,-c\right\}  $ and, so, $\operatorname*{dom}\ell
=\operatorname*{dom}f$ $(\neq\emptyset)\ $and the $\ell_{t}$'s are proper.
Since $f$ is lsc at $\theta$ (because it is subdifferentiable at $\theta$)
there exists an open\ neighborhood $U\in\mathcal{N}_{X}$ such that
\[
f(x)\geq-c\text{ \ \ \ for all }x\in U.
\]
Hence, $f=\ell$ and $\operatorname*{cl}\ell=\operatorname*{cl}f$ on $U$;
consequently, $\partial f(\theta)=\partial\ell(\theta)$ and we have proved (i)
and (iv). Now we proceed by proving (ii):
\begin{align*}
\operatorname*{cl}\ell &  =\operatorname*{cl}(\sup_{t\in T}(\max
\{f_{t},-c\}))\\
&  =\operatorname*{cl}(\max\{\sup_{t\in T}f_{t},-c\})\\
&  =\max\{\operatorname*{cl}(\sup_{t\in T}f_{t}),-c\}\text{ \ (by (\ref{cl}%
))}\\
&  =\max\{\sup_{t\in T}(\operatorname*{cl}f_{t}),-c\}\text{ \ (by
(\ref{assumption}))}\\
&  =\sup_{t\in T}\max\{\operatorname*{cl}f_{t},-c\}\\
&  =\sup_{t\in T}\operatorname*{cl}(\max\{f_{t},-c\})\text{ \ (by \ref{cl})}\\
&  =\sup_{t\in T}(\operatorname*{cl}\ell_{t}),
\end{align*}
leading to\ (ii).

Finally, to prove (iii), observe that for every $t\in T$ such that $\ell
_{t}(\theta)>-c$ we have
\[
\ell_{t}(\theta)=\max\{f_{t}(\theta),-c\}=f_{t}(\theta),
\]
and so for all $\varepsilon\in\left[  0,c\right[  $
\[
\{t\in T\mid\ell_{t}(\theta)\geq\ell(\theta)-\varepsilon\}=\{t\in T\mid
\ell_{t}(\theta)\geq-\varepsilon\}=\{t\in T\mid f_{t}(\theta)\geq
-\varepsilon\}=T_{\varepsilon}(\theta),
\]
yielding (iii).
\end{dem}

\begin{dem}
(\textbf{of\ Theorem \ref{thmcompact0}}) First, we show the inclusion
"$\supset$". We start by verifying that for every $t\in T$
\begin{equation}
\partial(f_{t}+\mathrm{I}_{\operatorname*{dom}f})(x)\subset\partial
(f_{t}+\mathrm{I}_{\overline{\operatorname*{dom}f}})(x).\label{ls}%
\end{equation}
We fix $x_{0}\in\operatorname*{ri}(\operatorname*{dom}f)$ and pick\ $z^{\ast
}\in\partial(f_{t}+\mathrm{I}_{\operatorname*{dom}f})(x).$ Given
$y\in\overline{\operatorname*{dom}f},$ we define
\[
y_{\lambda}:=\lambda x_{0}+(1-\lambda)y,\text{ \ }\lambda\in\left]
0,1\right[  ,
\]
so that $y_{\lambda}\in\operatorname*{dom}f$ by the accessibility lemma, and
\[
\left\langle z^{\ast},y_{\lambda}-x\right\rangle \leq f_{t}(y_{\lambda}%
)-f_{t}(x)\leq\lambda f(x_{0})+(1-\lambda)f_{t}(y)-f_{t}(x).
\]
As $\lambda\downarrow0$ we get
\begin{equation}
\left\langle z^{\ast},y-x\right\rangle \leq f_{t}(y)-f_{t}(x),\label{arg}%
\end{equation}
and, so, $z^{\ast}\in\partial(f_{t}+\mathrm{I}_{\overline{\operatorname*{dom}%
f}})(x).\medskip$

Next, for every $p\in\mathcal{P}$ and $\varepsilon>0$ we have, by
(\ref{ls})\ and (\ref{rel}),
\begin{align}
\overline{\operatorname*{co}}\left\{  \bigcup\limits_{t\in T(x)}\partial
(f_{t}+\mathrm{I}_{\operatorname*{dom}f})(x)\right\}   &  \subset
\overline{\operatorname*{co}}\left\{  \bigcup\limits_{t\in T(x)}\partial
(f_{t}+\mathrm{I}_{\overline{\operatorname*{dom}f}})(x)\right\} \nonumber\\
&  \subset\overline{\operatorname*{co}}\left\{  \bigcup\limits_{t\in
T_{\varepsilon}(x)}\breve{\partial}_{p}^{\varepsilon}(f_{t}+\mathrm{I}%
_{\overline{\operatorname*{dom}f}})(x)\right\}  . \label{inc}%
\end{align}
So, due to (\ref{fri}) we obtain
\[
\overline{\operatorname*{co}}\left\{  \bigcup\limits_{t\in T(x)}\partial
(f_{t}+\mathrm{I}_{\operatorname*{dom}f})(x)\right\}  \subset\bigcap
_{\varepsilon>0,\text{ }p\in\mathcal{P}}\overline{\operatorname*{co}}\left\{
\bigcup\limits_{t\in T_{\varepsilon}(x)}\breve{\partial}_{p}^{\varepsilon
}(f_{t}+\mathrm{I}_{\overline{\operatorname*{dom}f}})(x)\right\}  =\partial
f(x).
\]

To prove the inclusion "$\subset$", it suffices to consider the nontrivial
case $\partial f(x)\neq\emptyset,$ entaling that $x\in\operatorname*{dom}f$
and (see (\ref{noemptys}))%
\begin{equation}
f(x)=(\operatorname*{cl}f)(x)\text{ and }\partial_{\varepsilon}f(x)=\partial
_{\varepsilon}(\operatorname*{cl}f)(x)\text{ \ for all }\varepsilon\geq0;
\label{clp}%
\end{equation}
and because $f(x)\in\mathbb{R}$ we can suppose that $x=\theta$ and
$f(\theta)=(\operatorname*{cl}f)(\theta)=0.$ Moreover, as $\partial
f(\theta)\neq\emptyset,$ both functions $f$ and $\operatorname*{cl}f$ are
proper. In addition, since we have that $\operatorname*{ri}%
(\operatorname*{dom}f)=\operatorname*{ri}(\operatorname*{dom}%
(\operatorname*{cl}f)),$ $\operatorname*{aff}(\operatorname*{dom}%
f)=\operatorname*{aff}(\operatorname*{dom}(\operatorname*{cl}f))$, and
$\operatorname*{cl}f\leq f,$ we deduce that
\begin{equation}
(\operatorname*{cl}f)_{\mid\operatorname*{aff}(\operatorname*{dom}%
(\operatorname*{cl}f))}\text{ is finite and continuous on }\operatorname*{ri}%
(\operatorname*{dom}(\operatorname*{cl}f)). \label{ac}%
\end{equation}

In a first step, we suppose that $\operatorname*{cl}f_{t}$ is proper for all
$t\in T.$ We define\
\begin{equation}
g_{t}:=(\operatorname*{cl}f_{t})+\mathrm{I}_{\overline{\operatorname*{dom}f}%
},\text{ }t\in T,\label{gt}%
\end{equation}
so that, by (\ref{assumption}) and the relation $\operatorname*{dom}%
(\operatorname*{cl}f)\subset\overline{\operatorname*{dom}f}$,%
\begin{equation}
\sup_{t\in T}g_{t}=\sup_{t\in T}(\operatorname*{cl}f_{t})+\mathrm{I}%
_{\overline{\operatorname*{dom}f}}=\operatorname*{cl}f+\mathrm{I}%
_{\overline{\operatorname*{dom}f}}=\operatorname*{cl}f.\label{gf}%
\end{equation}
Thus, taking into account (\ref{ac}), we are in position to apply (\ref{fri})
to the lsc proper functions $g_{t},$ and we get
\begin{align}
\partial f(\theta) &  =\partial(\operatorname*{cl}f)(\theta)\text{ (by
(\ref{clp}))\smallskip}\nonumber\\
&  =\partial(\sup_{t\in T}g_{t})(\theta)\text{ (by (\ref{gf}))}\nonumber\\
&  =\bigcap_{_{\varepsilon>0,\text{ }p\in\mathcal{P}}}\overline
{\operatorname*{co}}\left\{  \bigcup\limits_{t\in\tilde{T}_{\varepsilon
}(\theta)}\breve{\partial}_{p}^{\varepsilon}(g_{t}+\mathrm{I}_{\overline
{\operatorname*{dom}(\operatorname*{cl}f)}})(\theta)\right\}  \text{ (by
(\ref{fri}))}\nonumber\\
&  =\bigcap_{_{\varepsilon>0,\text{ }p\in\mathcal{P}}}\overline
{\operatorname*{co}}\left\{  \bigcup\limits_{t\in\tilde{T}_{\varepsilon
}(\theta)}\breve{\partial}_{p}^{\varepsilon}g_{t}(\theta)\right\}  \text{
(since }\operatorname*{dom}g_{t}\subset\overline{\operatorname*{dom}%
f}=\overline{\operatorname*{dom}(\operatorname*{cl}f)}\text{)}\nonumber\\
&  =\bigcap_{_{\varepsilon>0,\text{ }p\in\mathcal{P}}}\overline
{\operatorname*{co}}\left\{  \bigcup\limits_{t\in\tilde{T}_{\varepsilon
}(\theta)}\breve{\partial}_{p}^{\varepsilon}g_{t}(\theta)\cap\partial
_{2\varepsilon}g_{t}(\theta)\right\}  \text{ (by (\ref{rel})),}\label{ar0a}%
\end{align}
where
\begin{align}
\tilde{T}_{\varepsilon}(\theta) &  :=\{t\in T\mid g_{t}(\theta)\geq(\sup_{t\in
T}g_{t})(\theta)-\varepsilon\}\nonumber\\
&  =\{t\in T\mid(\operatorname*{cl}f_{t})(\theta)\geq-\varepsilon\}\subset
T_{\varepsilon}(\theta).\label{nu1}%
\end{align}
Moreover, for every $t\in\tilde{T}_{\varepsilon}(\theta)$ we have that
$\partial_{2\varepsilon}g_{t}(\theta)\subset\partial_{3\varepsilon}f(\theta),$
which comes from the following inequalities: for $z\in\operatorname*{dom}f$
and $z^{\ast}\in\partial_{2\varepsilon}g_{t}(\theta),$%
\begin{align*}
f(z)\geq f_{t}(z)\geq(\operatorname*{cl}f_{t})(z)=g_{t}(z) &  \geq
g_{t}(\theta)+\left\langle z^{\ast},z\right\rangle -2\varepsilon\\
&  \geq(\sup_{t\in T}g_{t})(\theta)+\left\langle z^{\ast},z\right\rangle
-3\varepsilon\\
&  =(\operatorname*{cl}f)(\theta)+\left\langle z^{\ast},z\right\rangle
-3\varepsilon\text{ \ (by (\ref{gf}))}\\
&  =\left\langle z^{\ast},z\right\rangle -3\varepsilon.
\end{align*}
So, (\ref{ar0a}) yields
\begin{equation}
\partial f(\theta)=\bigcap_{_{\varepsilon>0,\text{ }p\in\mathcal{P}}}%
\overline{\operatorname*{co}}\left\{  \bigcup\limits_{t\in\tilde
{T}_{\varepsilon}(\theta)}\breve{\partial}_{p}^{\varepsilon}g_{t}(\theta
)\cap\partial_{3\varepsilon}f(\theta)\right\}  .\label{ar0}%
\end{equation}

Take\ $x^{\ast}\in\partial f(\theta)$ and fix $u\in X.$ Given $k\in\mathbb{N}$
and $p\in\mathcal{P}$, (\ref{ar0}) ensures the following relation in
$\mathbb{R}^{2}$
\begin{align*}
\left(
\begin{array}
[c]{c}%
\left\langle x^{\ast},u\right\rangle \\
\left\langle x^{\ast},x_{0}\right\rangle
\end{array}
\right)   &  \in\overline{\operatorname*{co}}\left\{  \bigcup\limits_{t\in
\tilde{T}_{\frac{1}{k}}(\theta)}\left\{  \left(
\begin{array}
[c]{c}%
\left\langle y^{\ast},u\right\rangle \\
\left\langle y^{\ast},x_{0}\right\rangle
\end{array}
\right)  \left\vert ~y^{\ast}\in\breve{\partial}_{p}^{\frac{1}{k}}g_{t}%
(\theta)\cap\partial_{\frac{3}{k}}f(\theta)\right.  \right\}  \right\}  \\
&  \subset\operatorname*{co}\left\{  \bigcup\limits_{t\in\tilde{T}_{\frac
{1}{k}}(\theta)}\left\{  \left(
\begin{array}
[c]{c}%
\left\langle y^{\ast},u\right\rangle \\
\left\langle y^{\ast},x_{0}\right\rangle
\end{array}
\right)  \left\vert ~y^{\ast}\in\breve{\partial}_{p}^{\frac{1}{k}}g_{t}%
(\theta)\cap\partial_{\frac{3}{k}}f(\theta)\right.  \right\}  \right\}
+\left[  -\frac{1}{k},\frac{1}{k}\right]  ^{2}.
\end{align*}
Therefore, taking into account Charathéodory's Theorem,\ there are%
\[
(\lambda_{1,k},\lambda_{2,k},\lambda_{3,k})\in\Delta_{3}:=\{(\lambda
_{1},\lambda_{2},\lambda_{3})\mid\lambda_{i}\geq0,\lambda_{1}+\lambda
_{2}+\lambda_{3}=1\}\text{ }%
\]
and
\begin{equation}
t_{i,k}\in\tilde{T}_{1/k}(\theta),\text{ and }y_{i,k}^{\ast}\in\breve
{\partial}_{p}^{1/k}g_{t_{i,k}}(\theta)\cap\partial_{3/k}f(\theta)\text{,
}i=1,2,3,\label{hya}%
\end{equation}
such that
\begin{equation}
\left\langle x^{\ast},u\right\rangle \leq\left\langle \lambda_{1,k}%
y_{1,k}^{\ast}+\lambda_{2,k}y_{2,k}^{\ast}+\lambda_{3,k}y_{3,k}^{\ast
},u\right\rangle +1/k,\label{be0}%
\end{equation}%
\begin{equation}
\left\vert \left\langle x^{\ast},x_{0}\right\rangle -\left\langle
\lambda_{1,k}y_{1,k}^{\ast}+\lambda_{2,k}y_{2,k}^{\ast}+\lambda_{3,k}%
y_{3,k}^{\ast},x_{0}\right\rangle \right\vert \leq1/k.\label{be}%
\end{equation}
Then, by the definition of $\breve{\partial}_{p}^{1/k},$ there exist
$y_{i,k}\in\operatorname*{dom}g_{t_{i,k}}\subset\overline{\operatorname*{dom}%
f}$ $(i=1,2,3)$ such that
\begin{equation}
p(y_{i,k})\leq1/k,\text{ }\left\vert (\operatorname*{cl}f_{t_{i,k}}%
)(y_{i,k})-(\operatorname*{cl}f_{t_{i,k}})(\theta)\right\vert \leq1/k,\text{
}\left\vert \left\langle y_{i,k}^{\ast},y_{i,k}\right\rangle \right\vert
\leq1/k,\label{hy0}%
\end{equation}
and (recall (\ref{hya}))
\begin{equation}
y_{i,k}^{\ast}\in\partial g_{t_{i,k}}(y_{i,k})\cap\partial_{3/k}%
f(\theta).\label{hy}%
\end{equation}
From the continuity assumption of $f_{\mid\operatorname*{aff}%
(\operatorname*{dom}f)}$ at $x_{0}\in\operatorname*{ri}(\operatorname*{dom}f)$
we choose\ $m\geq0$ and\ $W\in\mathcal{N}_{X}$ such that\
\[
x_{0}+W\cap\operatorname*{aff}(\operatorname*{dom}f)\subset\operatorname*{dom}%
f\text{ and }\sup_{w\in W\cap\operatorname*{aff}(\operatorname*{dom}f)}%
f(x_{0}+w)\leq m.
\]
Then, using (\ref{hy}), for each $i=1,2,3$ and for all $z\in W\cap
\operatorname*{aff}(\operatorname*{dom}f)$
\begin{equation}
\langle y_{i,k}^{\ast},x_{0}+z\rangle\leq f(x_{0}+z)-f(\theta)+3/k=f(x_{0}%
+z)+3/k\leq m+3/k,\label{cc}%
\end{equation}
and so, since $\theta\in W\cap\operatorname*{aff}(\operatorname*{dom}f)$,
\begin{equation}
\langle y_{i,k}^{\ast},x_{0}\rangle\leq m+3/k.\label{bb}%
\end{equation}
Thus, from (\ref{be}) we deduce that the sequences $(\langle\lambda
_{i,k}y_{i,k}^{\ast},x_{0}\rangle)_{k},$ $i=1,2,3,$ are\ bounded.
Consequently, by (\ref{cc}) and (\ref{bb}) we see that there exists a positive
number $r$ such that
\begin{equation}
\langle\lambda_{i,k}y_{i,k}^{\ast},z\rangle\leq\lambda_{i,k}(m+3/k-\langle
y_{i,k}^{\ast},x_{0}\rangle)\leq r\text{ \ for all }z\in W\cap
\operatorname*{aff}(\operatorname*{dom}f);\label{dd}%
\end{equation}
that is,
\[
r^{-1}(\lambda_{i,k}y_{i,k}^{\ast})_{k}\subset(W\cap\operatorname*{aff}%
(\operatorname*{dom}f))^{\circ}:=\left\{  u^{\ast}\in X^{\ast}\mid\langle
u^{\ast},z\rangle\leq1\text{ }\forall z\in W\cap\operatorname*{aff}%
(\operatorname*{dom}f)\right\}  .
\]

We endow the lcs space $Y:=\operatorname*{aff}(\operatorname*{dom}f)$
$(=\overline{\operatorname*{aff}}(\operatorname*{dom}f))$ with the induced
topology from\ $X,$\ and denote by $Y^{\ast}$ its topological dual space. By
the Alaoglu-Bourbaki\ Theorem applied to the dual pair $(Y,Y^{\ast})$, there
exists a subnet of the sequence of the restrictions of $(\lambda_{i,k}%
y_{i,k}^{\ast})_{k}$ to $Y,$ denoted by $(\lambda_{i,k_{\alpha}}%
y_{i,k_{\alpha}\mid Y}^{\ast})_{\alpha\in\Upsilon},$ which weak*-converges to
some $\hat{y}_{i}^{\ast}\in Y^{\ast}$ $(i=1,2,3)$. Thus, if $y_{i}^{\ast}\in
X^{\ast}$ is an extension of $\hat{y}_{i}^{\ast}$ to $X,$ it satisfies for
every $z\in\overline{\operatorname*{dom}f}$
\begin{align}
\left\langle y_{i}^{\ast},z\right\rangle  &  =\left\langle \hat{y}_{i}^{\ast
},z\right\rangle =\lim_{\alpha\in\Upsilon}\left\langle \lambda_{i,k_{\alpha}%
}y_{i,k_{\alpha}}^{\ast},z\right\rangle \nonumber\\
&  =\lim_{\alpha\in\Upsilon}(\left\langle \lambda_{i,k_{\alpha}}%
y_{i,k_{\alpha}}^{\ast},z-y_{i,k_{\alpha}}\right\rangle +\left\langle
\lambda_{i,k_{\alpha}}y_{i,k_{\alpha}}^{\ast},y_{i,k_{\alpha}}\right\rangle
)\nonumber\\
&  \leq\liminf_{\alpha\in\Upsilon}\left(  \lambda_{i,k_{\alpha}}%
(g_{t_{i,k_{\alpha}}}(z)-g_{t_{i,k_{\alpha}}}(y_{i,k_{\alpha}}))+1/k_{\alpha
}\right)  \text{ \ (by (\ref{hy}) and (\ref{hy0}))}\nonumber\\
&  =\liminf_{\alpha\in\Upsilon}\left(  \lambda_{i,k_{\alpha}}%
((\operatorname*{cl}f_{t_{i,k_{\alpha}}})(z)-(\operatorname*{cl}%
f_{t_{i,k_{\alpha}}})(y_{i,k_{\alpha}}))+1/k_{\alpha}\right) \nonumber\\
&  \leq\liminf_{\alpha\in\Upsilon}\left(  \lambda_{i,k_{\alpha}}%
((\operatorname*{cl}f_{t_{i,k_{\alpha}}})(z)-(\operatorname*{cl}%
f_{t_{i,k_{\alpha}}})(\theta))+2/k_{\alpha}\right)  \text{ \ (by (\ref{hy0}%
))}\nonumber\\
&  \leq\liminf_{\alpha\in\Upsilon}\left(  \lambda_{i,k_{\alpha}}%
(\operatorname*{cl}f_{t_{i,k_{\alpha}}})(z)+3/k_{\alpha}\right)  \text{ \ (as
}t_{i,k_{\alpha}}\in\tilde{T}_{1/k_{\alpha}}(\theta),\text{ by (\ref{hya}%
))}\nonumber\\
&  =\liminf_{\alpha\in\Upsilon}\lambda_{i,k_{\alpha}}(\operatorname*{cl}%
f_{t_{i,k_{\alpha}}})(z)\label{ee0}\\
&  \leq\liminf_{\alpha\in\Upsilon}\lambda_{i,k_{\alpha}}f_{t_{i,k_{\alpha}}%
}(z), \label{ee}%
\end{align}
where $(k_{\alpha})_{\alpha\in\Upsilon}\subset\mathbb{N}$ is such that
$\lim_{\alpha\in\Upsilon}k_{\alpha}=+\infty$ and $1/k_{\alpha}\leq
\varepsilon_{0},$ eventually.

We may suppose that $\lambda_{i,k_{\alpha}}\rightarrow\lambda_{i}\ $for some
$\lambda_{i}\in\left[  0,1\right]  ,$ so that $\lambda:=(\lambda_{1}%
,\lambda_{2},\lambda_{3})\in\Delta_{3}.$ We denote $I_{0}:=\{i=1,2,3\mid
\lambda_{i}=0\}$ (this set can be empty). If $i\in I_{0},$ from (\ref{ee}) we
get
\[
\left\langle y_{i}^{\ast},z\right\rangle \leq\liminf_{\alpha\in\Upsilon
}\lambda_{i,k_{\alpha}}f_{t_{i,k_{\alpha}}}(z)\leq\liminf_{\alpha\in\Upsilon
}\lambda_{i,k_{\alpha}}f(z)=0\text{ \ for all }z\in\operatorname*{dom}f,
\]
showing that
\begin{equation}
y_{i}^{\ast}\in\mathrm{N}_{\operatorname*{dom}f}(\theta)=\mathrm{N}%
_{\overline{\operatorname*{dom}f}}(\theta).\label{w1}%
\end{equation}
Otherwise, if $i\not \in I_{0}$ ($\left\{  1,2,3\right\}  \setminus I_{0}$ is
nonempty because $\sum_{i}\lambda_{i}=1$), (\ref{ee0}) gives rise to
\begin{equation}
\left\langle \lambda_{i}^{-1}y_{i}^{\ast},z\right\rangle \leq\liminf
_{\alpha\in\Upsilon}(\operatorname*{cl}f_{t_{i,k_{\alpha}}})(z)\leq
\liminf_{\alpha\in\Upsilon}f_{t_{i,k_{\alpha}}}(z)\text{ \ for all }%
z\in\overline{\operatorname*{dom}f}.\label{ff}%
\end{equation}
But, since $1/k_{\alpha}\leq\varepsilon_{0}$ eventually, by (\ref{nu1})
\begin{equation}
t_{i,k_{\alpha}}\in\tilde{T}_{1/k_{\alpha}}(\theta)\subset T_{1/k_{\alpha}%
}(\theta)\subset T_{\varepsilon_{0}}(\theta),\label{com}%
\end{equation}
the last\ set being compact by assumption (i), we may assume that
$(t_{i,k_{\alpha}})_{\alpha\in\Upsilon}$ converges to some $t_{i}\in
T_{\varepsilon_{0}}(\theta).$ Thus,
\begin{align*}
f(\theta)\geq f_{t_{i}}(\theta) &  \geq\limsup_{\alpha\in\Upsilon
}f_{t_{i,k_{\alpha}}}(\theta)\text{ \ (by (ii))}\\
&  \geq\limsup_{\alpha\in\Upsilon}(\operatorname*{cl}f_{t_{i,k_{\alpha}}%
})(\theta)\\
&  \geq\limsup_{\alpha\in\Upsilon}(f(\theta)-1/k_{\alpha})=f(\theta)\text{
\ (by (\ref{com}))};
\end{align*}
that is, $t_{i}\in T(\theta).$ So, (\ref{ff}) together with hypothesis (ii)
yield for all $z\in\operatorname*{dom}f$
\[
\left\langle \lambda_{i}^{-1}y_{i}^{\ast},z\right\rangle \leq\liminf
_{\alpha\in\Upsilon}f_{t_{i,k_{\alpha}}}(z)\leq\limsup_{\alpha\in\Upsilon
}f_{t_{i,k_{\alpha}}}(z)\leq f_{t_{i}}(z)=f_{t_{i}}(z)-f_{t_{i}}(\theta);
\]
so that
\begin{equation}
\lambda_{i}^{-1}y_{i}^{\ast}\in\partial(f_{t_{i}}+\mathrm{I}%
_{\operatorname*{dom}f})(\theta).\label{w2}%
\end{equation}
By taking limits in (\ref{be0}) when $u\in Y$ we get (recall that $k_{\alpha
}\rightarrow+\infty$)
\begin{align*}
\left\langle x^{\ast},u\right\rangle  &  \leq\limsup_{\alpha\in\Upsilon
}\left(  \left\langle \lambda_{1,k_{\alpha}}y_{1,k_{\alpha}}^{\ast}%
+\lambda_{2,k_{\alpha}}y_{2,k_{\alpha}}^{\ast}+\lambda_{3,k_{\alpha}%
}y_{3,k_{\alpha}}^{\ast},u\right\rangle +1/k_{\alpha}\right)  \\
&  \leq\sum_{i\in I_{0}}\limsup_{\alpha\in\Upsilon}\left\langle \lambda
_{i,k_{\alpha}}y_{i,k_{\alpha}}^{\ast},u\right\rangle +\sum_{i\not \in I_{0}%
}\limsup_{\alpha\in\Upsilon}\left\langle \lambda_{i,k_{\alpha}}y_{i,k_{\alpha
}}^{\ast},u\right\rangle \\
&  =\sum_{i\in I_{0}}\left\langle y_{i}^{\ast},u\right\rangle +\left\langle
\sum_{i\not \in I_{0}}\lambda_{i}(\lambda_{i}^{-1}y_{i}^{\ast}),u\right\rangle
\text{ (as }(\lambda_{i,k_{\alpha}}y_{i,k_{\alpha}}^{\ast})_{\mid Y}%
\overset{w^{\ast}}{\longrightarrow}\hat{y}_{i}^{\ast}=(y_{i}^{\ast})_{\mid
Y}\text{),}%
\end{align*}
which gives us, due to (\ref{w1}) and (\ref{w2}),
\begin{equation}
\left\langle x^{\ast},u\right\rangle \leq\mathrm{\sigma}_{\mathrm{N}%
_{\operatorname*{dom}f}(\theta)}(u)+\mathrm{\sigma}_{\sum_{i\not \in I_{0}%
}\lambda_{i}\partial(f_{t_{i}}+\mathrm{I}_{\operatorname*{dom}f})(\theta
)}(u).\label{he}%
\end{equation}
Observe that when $u\not \in Y$ $(=\operatorname*{aff}(\operatorname*{dom}%
f)),$ we have that $u\not \in \overline{\mathbb{R}_{+}(\operatorname*{dom}%
f)}=(\mathrm{N}_{\operatorname*{dom}f}(\theta))^{\bot}$ and so
\[
\mathrm{\sigma}_{\mathrm{N}_{\operatorname*{dom}f}(\theta)}(u)=+\infty;
\]
hence, thanks to (\ref{he}), for every $u\in X$ it holds
\begin{align}
\left\langle x^{\ast},u\right\rangle  &  \leq\mathrm{\sigma}_{\mathrm{N}%
_{\operatorname*{dom}f}(\theta)}(u)+\mathrm{\sigma}_{\sum_{i\not \in I_{0}%
}\lambda_{i}\partial(f_{t_{i}}+\mathrm{I}_{\operatorname*{dom}f})(\theta
)}(u)\nonumber\\
&  =\mathrm{\sigma}_{\mathrm{N}_{\operatorname*{dom}f}(\theta)+\sum
_{i\not \in I_{0}}\lambda_{i}\partial(f_{t_{i}}+\mathrm{I}%
_{\operatorname*{dom}f})(\theta)}(u).\label{rr}%
\end{align}
Additionally,\ for\ each $i\not \in I_{0}$
\[
\mathrm{N}_{\operatorname*{dom}f}(\theta)+\lambda_{i}\partial(f_{t_{i}%
}+\mathrm{I}_{\operatorname*{dom}f})(\theta)\subset\partial(\lambda
_{i}f_{t_{i}}+\mathrm{I}_{\operatorname*{dom}f})(\theta)=\lambda_{i}%
\partial(f_{t_{i}}+\mathrm{I}_{\operatorname*{dom}f})(\theta),
\]
and so (\ref{rr}) yields (recall that $t_{i}\in T(\theta)$)%
\[
\left\langle x^{\ast},u\right\rangle \leq\mathrm{\sigma}_{\sum_{i\in
\{1,2,3\}\setminus I_{0}}\lambda_{i}\partial(f_{t_{i}}+\mathrm{I}%
_{\operatorname*{dom}f})(\theta)}(u)\leq\mathrm{\sigma}_{\overline
{\operatorname*{co}}\left\{  \bigcup\nolimits_{t\in T(\theta)}\partial
(f_{t}+\mathrm{I}_{\operatorname*{dom}f})(\theta)\right\}  }(u).
\]
As $u$ is arbitrary in $X$ we conclude that
\[
x^{\ast}\in\overline{\operatorname*{co}}\left\{  \bigcup\limits_{t\in
T(\theta)}\partial(f_{t}+\mathrm{I}_{\operatorname*{dom}f})(\theta)\right\}  ,
\]
proving the desired inclusion when the $\operatorname*{cl}f_{t}$'s are proper.

Finally, to deal with the case when the $(\operatorname*{cl}f_{t})$'s are not
all necessarily proper, we consider the functions%
\[
\ell_{t}:=\max\{f_{t},-2\varepsilon_{0}\},\text{ }t\in T,\text{ and }%
\ell:=\sup_{t\in T}\ell_{t}.
\]
According to Lemma \ref{lt}, the $\operatorname*{cl}\ell_{t}$'s are proper,
$\left\{  t\in T\mid\ell_{t}(\theta)\geq\ell(\theta)-\varepsilon\right\}
=T_{\varepsilon}(\theta)$ for all $\varepsilon\in\left[  0,\varepsilon
_{0}\right[  $ $(\subset\left[  0,2\varepsilon_{0}\right[  ),$
$\operatorname*{dom}\ell=\operatorname*{dom}f,$ $\operatorname*{cl}\ell
:=\sup_{t\in T}(\operatorname*{cl}\ell_{t}),$ and $\partial f(\theta
)=\partial\ell(\theta).$ It follows that the family $\left\{  \ell_{t},\text{
}t\in T\right\}  $ satisfies the current conditions (i) and (ii).
Consequently, from the previous part of the proof, applied to the functions
$\ell_{t},$ $t\in T,$ we get (observe that $\ell(\theta)=\max\{f_{t}%
(\theta),-2\varepsilon_{0}\}=f(\theta)=0$)%
\[
\partial f(\theta)=\partial\ell(\theta)=\overline{\operatorname*{co}}\left\{
\bigcup\limits_{t\in\left\{  t\in T\mid\ell_{t}(\theta)=0\right\}  }%
\partial(\ell_{t}+\mathrm{I}_{\operatorname*{dom}\ell})(\theta)\right\}
=\overline{\operatorname*{co}}\left\{  \bigcup\limits_{t\in T(\theta)}%
\partial(\ell_{t}+\mathrm{I}_{\operatorname*{dom}f})(\theta)\right\}  .
\]
Take $t\in T(\theta)$ such that $\partial(\ell_{t}+\mathrm{I}%
_{\operatorname*{dom}\ell})(\theta)\neq\emptyset$ (such a $t$ always exists
because $\partial f(\theta)\neq\emptyset$). We show that
\begin{equation}
\partial(\ell_{t}+\mathrm{I}_{\operatorname*{dom}f})(\theta)=\partial
(f_{t}+\mathrm{I}_{\operatorname*{dom}f})(\theta). \label{fi}%
\end{equation}
Since $\ell_{t}+\mathrm{I}_{\operatorname*{dom}\ell}=\max\{f_{t}%
+\mathrm{I}_{\operatorname*{dom}\ell},-2\varepsilon_{0}\}$ and $(f_{t}%
+\mathrm{I}_{\operatorname*{dom}\ell})(\theta)=0>-2\varepsilon_{0},$ it
suffices to verify that $f_{t}+\mathrm{I}_{\operatorname*{dom}\ell}$ is lsc at
$\theta,$ because then the two functions $\ell_{t}+\mathrm{I}%
_{\operatorname*{dom}\ell}$ and $f_{t}+\mathrm{I}_{\operatorname*{dom}\ell}$
coincide in a neighborhood of $\theta.$ Indeed, using (\ref{cl}) and the lower
semicontinuity of $\ell_{t}+\mathrm{I}_{\operatorname*{dom}\ell}$ at $\theta$
(a consequence of the nonemptiness of $\partial(\ell_{t}+\mathrm{I}%
_{\operatorname*{dom}\ell})(\theta)$), we have
\begin{align*}
\max\{\operatorname*{cl}(f_{t}+\mathrm{I}_{\operatorname*{dom}\ell}%
)(\theta),-2\varepsilon_{0}\}  &  =\operatorname*{cl}(\max\{f_{t}%
+\mathrm{I}_{\operatorname*{dom}\ell},-2\varepsilon_{0}\})(\theta)\\
&  =\operatorname*{cl}(\ell_{t}+\mathrm{I}_{\operatorname*{dom}\ell}%
)(\theta)=(\ell_{t}+\mathrm{I}_{\operatorname*{dom}\ell})(\theta)=0;
\end{align*}
that is, $\operatorname*{cl}(f_{t}+\mathrm{I}_{\operatorname*{dom}\ell
})(\theta)=0.$ But $t\in T(\theta),$ and so $(f_{t}+\mathrm{I}%
_{\operatorname*{dom}\ell})(\theta)=0=\operatorname*{cl}(f_{t}+\mathrm{I}%
_{\operatorname*{dom}\ell})(\theta),$ yielding the lower semicontinuity of
$\ell_{t}+\mathrm{I}_{\operatorname*{dom}\ell}$ at $\theta,$ and (\ref{fi}) follows.

The proof of the inclusion "$\subset$" is done.
\end{dem}

\begin{rem}
\emph{In Theorem \ref{thmcompact0} one can replace (i)-(ii) by\ the following
weaker pair of conditions (i}$^{\prime}$\emph{)-(ii}$^{\prime}$\emph{), which
emphasizes the role played by the }$\varepsilon$\emph{-active sets at }$x:$

$(i^{\prime})$\emph{\ the sets} $T_{\varepsilon}(x)$ \emph{are compact for
all} $\varepsilon\in\left[  0,\varepsilon_{0}\right]  ,$

$(ii^{\prime})$\emph{\ the functions} $t\mapsto f_{t}(z),$ $z\in
\operatorname*{dom}f,$\ \emph{are usc on} $T(x).$ \newline\emph{Indeed,
from\ (\ref{com}) we have that} $t_{i,k_{\alpha}}\in\tilde{T}_{1/k_{\alpha}%
}(x)\subset T_{1/k_{\alpha}}(x)$ \emph{for all} $\alpha\in\Upsilon,$ \emph{and
so }$(i^{\prime})$\emph{ gives rise to }%
\[
t_{i}\in\cap_{\alpha}\operatorname*{cl}(T_{1/k_{\alpha}}(x))=\cap_{\alpha
}T_{1/k_{\alpha}}(x)=T(x).
\]
\emph{Thus, the proof follows by using the weaker assumption}\ $(ii^{\prime
}),$ \emph{instead of} $(ii).$
\end{rem}

\begin{rem}
\label{remcl}\emph{\cite[Corollary 9]{HanLop08}} \emph{The closure condition
(\ref{assumption}) holds in each one of the following situations:}

\emph{(1) the functions }$f_{t}$\emph{, }$t\in T,$\emph{ are lsc. }

\emph{(2) the }$f_{t}$\emph{'s have a common continuity point; this follows
if, for instance, the supremum function }$f$\emph{ is finite and continuous at
some point.}

\emph{(3) }$T$\emph{ is finite and all but one of the functions }$f_{t}%
$\emph{'s have a common continuity point in }$\operatorname*{dom}f$
\emph{(this includes (\ref{cl}))}$.$

\emph{(4) }$X$\emph{ is finite-dimensional, and the relative interiors
}$\operatorname*{ri}(\operatorname*{dom}f_{t}),$\emph{ }$t\in T,$\emph{ have a
common point in\ }$\operatorname*{dom}f.$
\end{rem}

\begin{theo}
\label{thmcompact1}Assume that the convex functions $f_{t},$ $t\in T,$ are
proper and lsc. Let $x\in X$ be such that for some $\varepsilon_{0}>0$\emph{:}

$(i)$ the set $T_{\varepsilon_{0}}(x)$ is compact,

$(ii)$ the functions $t\mapsto f_{t}(z),$ $z\in\operatorname*{dom}f,$\ are usc
on $T_{\varepsilon_{0}}(x).$ \newline Then
\[
\partial f(x)=\bigcap_{L\in\mathcal{F}(x)}\overline{\operatorname*{co}%
}\left\{  \bigcup\limits_{t\in T(x)}\partial(f_{t}+\mathrm{I}_{L\cap
\operatorname*{dom}f})(x)\right\}  ,
\]
where $\mathcal{F}(x):=\{$finite-dimensional linear subspaces of $X$ such that
$x\in L\}$.
\end{theo}

\begin{dem}
We start by verifying the inclusion "$\supset$". Fix $L\in\mathcal{F}(x)$,
$p\in\mathcal{P}$ and $\varepsilon>0.$ By arguing as in the proof of the
inclusion "$\supset$" in Theorem \ref{thmcompact0} (see (\ref{inc})) we can
show\ that
\begin{align*}
\overline{\operatorname*{co}}\left\{  \bigcup\limits_{t\in T(x)}\partial
(f_{t}+\mathrm{I}_{L\cap\operatorname*{dom}f})(x)\right\}   &  \subset
\overline{\operatorname*{co}}\left\{  \bigcup\limits_{t\in T(x)}\partial
(f_{t}+\mathrm{I}_{\overline{L\cap\operatorname*{dom}f}})(x)\right\}  \\
&  \subset\overline{\operatorname*{co}}\left\{  \bigcup\limits_{t\in
T_{\varepsilon}(x)}\breve{\partial}_{p}^{\varepsilon}(f_{t}+\mathrm{I}%
_{\overline{L\cap\operatorname*{dom}f}})(x)\right\}  .
\end{align*}
So, due to (\ref{mn}) we obtain
\begin{align*}
\bigcap_{L\in\mathcal{F}(x)}\overline{\operatorname*{co}}\left\{
\bigcup\limits_{t\in T(x)}\partial(f_{t}+\mathrm{I}_{L\cap\operatorname*{dom}%
f})(x)\right\}   &  \subset\bigcap_{\substack{\varepsilon>0,\text{ }%
p\in\mathcal{P}\\L\in\mathcal{F}(x)}}\overline{\operatorname*{co}}\left\{
\bigcup\limits_{t\in T_{\varepsilon}(x)}\breve{\partial}_{p}^{\varepsilon
}(f_{t}+\mathrm{I}_{\overline{L\cap\operatorname*{dom}f}})(x)\right\}
\smallskip\\
&  =\partial f(x).
\end{align*}
Observe that this argument does not use the lower semicontinuity of the
$f_{t}$'s.

We are going to prove the inclusion "$\subset$" by considering the nontrivial
case $\partial f(x)\neq\emptyset,$ entaling that $x\in\operatorname*{dom}f$.
Let us suppose that $x=\theta$, $f(\theta)=(\operatorname*{cl}f)(\theta)=0.$

For a fixed $L\in\mathcal{F}(\theta)$ we apply Theorem \ref{thmcompact0}\ to
the lsc proper convex functions
\[
h_{t}:=f_{t}+\mathrm{I}_{\overline{L\cap\operatorname*{dom}f}},\text{ }t\in T,
\]%
\[
h:=\sup_{t\in T}h_{t}=f+\mathrm{I}_{\overline{L\cap\operatorname*{dom}f}}.
\]
Obviously, $h\geq f,$ $\operatorname*{dom}h\subset\operatorname*{dom}f,$
$h(\theta)=f(\theta)=0,$ $\operatorname*{ri}(\operatorname*{dom}%
h)\neq\emptyset$ (because $\operatorname*{dom}h$ $\subset L$), and so we have
that $h_{\mid\operatorname*{aff}(\operatorname*{dom}h)}$ is continuous on
$\operatorname*{ri}(\operatorname*{dom}h),$ together with
\begin{equation}
\partial f(\theta)\subset\partial h(\theta)\label{num}%
\end{equation}
and
\begin{equation}
\{t\in T\mid h_{t}(\theta)\geq h(\theta)-\varepsilon\}=\{t\in T\mid
f_{t}(\theta)\geq-\varepsilon\}=T_{\varepsilon}(\theta)\text{ \ for all
}\varepsilon\geq0.\label{te}%
\end{equation}
Since, for $z\in\operatorname*{dom}h$ $(\subset\operatorname*{dom}f)$ the
function $t\mapsto h_{t}(z)=f_{t}(z)+\mathrm{I}_{\overline{L\cap
\operatorname*{dom}f}}(z)$ is usc on $T_{\varepsilon_{0}}(\theta),$ Theorem
\ref{thmcompact0} applies and yields, taking into account (\ref{te}),\
\begin{align*}
\partial h(\theta) &  =\overline{\operatorname*{co}}\left\{  \bigcup
\limits_{\{t\in T\mid h_{t}(\theta)=h(\theta)\}}\partial(h_{t}+\mathrm{I}%
_{\operatorname*{dom}h})(\theta)\right\}  \\
&  =\overline{\operatorname*{co}}\left\{  \bigcup\limits_{t\in T(\theta
)}\partial(f_{t}+\mathrm{I}_{\overline{L\cap\operatorname*{dom}f}}%
+\mathrm{I}_{\operatorname*{dom}f\cap\overline{L\cap\operatorname*{dom}f}%
})(\theta)\right\}  \\
&  =\overline{\operatorname*{co}}\left\{  \bigcup\limits_{t\in T(\theta
)}\partial(f_{t}+\mathrm{I}_{\overline{L\cap\operatorname*{dom}f}}%
+\mathrm{I}_{L\cap\operatorname*{dom}f\cap\overline{L\cap\operatorname*{dom}%
f}})(\theta)\right\}  ,
\end{align*}
and hence, by (\ref{num}),
\begin{equation}
\partial f(\theta)\subset\partial h(\theta)\subset\overline{\operatorname*{co}%
}\left\{  \bigcup\limits_{t\in T(\theta)}\partial(f_{t}+\mathrm{I}%
_{L\cap\operatorname*{dom}f})(\theta)\right\}  .\label{pp}%
\end{equation}
Thus, due to\ the arbitrariness of the $L$'s, the desired inclusion follows.
\end{dem}

\begin{theo}
\label{thmcompact}Assume that the convex functions $f_{t},$ $t\in T,$ satisfy
\[
\operatorname*{cl}f=\sup_{t\in T}(\operatorname*{cl}f_{t}).
\]
Let $x\in X$ be such that for some $\varepsilon_{0}>0$\emph{:}

$(i)$ the set $T_{\varepsilon_{0}}(x)$ is compact,

$(ii)$ the functions $t\mapsto\left(  \operatorname*{cl}f_{t}\right)  (z),$
$z\in\operatorname*{dom}\left(  \operatorname*{cl}f\right)  ,$\ are usc on
$T_{\varepsilon_{0}}(x).$ \newline Then
\begin{equation}
\partial f(x)=\bigcap_{L\in\mathcal{F}(x)}\overline{\operatorname*{co}%
}\left\{  \bigcup\limits_{t\in T(x)}\partial(f_{t}+\mathrm{I}_{L\cap
\operatorname*{dom}f})(x)\right\}  . \label{form5}%
\end{equation}

\end{theo}

\begin{dem}
The inclusion "$\supset$" is proved as\ in\ Theorem \ref{thmcompact1}. To
prove the inclusion "$\subset$" we suppose that $\partial f(x)\neq\emptyset,$
allowing us to take $x=\theta$ and $f(\theta)=0.$ Let us\ consider the proper
convex functions
\[
\ell_{t}:=\max\{f_{t},-2\varepsilon_{0}\},\text{ }t\in T,
\]
and
\begin{equation}
\ell:=\sup_{t\in T}\ell_{t}=\max\left\{  f,-2\varepsilon_{0}\right\}
.\label{ppo}%
\end{equation}
According to Lemma \ref{lt}, in some open neighborhood $U$ of $\theta$ we have
that $f\geq-2\varepsilon_{0}$ and $f$ coincides with $\ell$, entailing
$(\emptyset\neq)$ $\partial f(\theta)=\partial\ell(\theta)\ $and
$(\operatorname*{cl}\ell)(\theta)=\ell(\theta)=f(\theta)=0.$ In addition, we
have $\operatorname*{dom}\ell=\operatorname*{dom}f$ and
\begin{equation}
\operatorname*{cl}\ell=\operatorname*{cl}f\text{ \ on }U.\label{dl}%
\end{equation}
From (\ref{ppo}) we also have that $\operatorname*{cl}\ell=\max\left\{
\operatorname*{cl}f,-2\varepsilon_{0}\right\}  $, and so $\operatorname*{dom}%
(\operatorname*{cl}\ell)=\operatorname*{dom}(\operatorname*{cl}f).$

Now, taking into account that $\operatorname*{cl}\ell=\sup_{t\in
T}(\operatorname*{cl}\ell_{t})$ by Lemma \ref{lt}(ii), we apply Theorem
\ref{thmcompact1} to the family $\{\operatorname*{cl}\ell_{t},$ $t\in T\}.$ To
this aim we need to verify the conditions of that theorem. Indeed, it is clear
that each $\operatorname*{cl}\ell_{t}$ is a proper convex lsc function.
Moreover, we have that for all $\varepsilon\in\left[  0,\varepsilon
_{0}\right]  $
\begin{align}
T_{\varepsilon}^{\prime}(\theta) &  :=\left\{  t\in T\mid(\operatorname*{cl}%
\ell_{t})(\theta)\geq(\operatorname*{cl}\ell)(\theta)-\varepsilon\right\}
\nonumber\\
&  =\left\{  t\in T\mid(\operatorname*{cl}\ell_{t})(\theta)=\operatorname*{cl}%
(\max\left\{  f_{t}(\theta),-2\varepsilon_{0}\right\}  )=\max\left\{
(\operatorname*{cl}f_{t})(\theta),-2\varepsilon_{0}\right\}  \geq
-\varepsilon\right\}  \nonumber\\
&  =\left\{  t\in T\mid(\operatorname*{cl}f_{t})(\theta)\geq-\varepsilon
\right\}  \subset T_{\varepsilon}(\theta).\label{di}%
\end{align}
Then, since that for every $z\in\operatorname*{dom}(\operatorname*{cl}%
\ell)=\operatorname*{dom}(\operatorname*{cl}f)$ the function $t\mapsto
(\operatorname*{cl}\ell_{t})(z)$ is usc on $T_{\varepsilon_{0}}^{\prime
}(\theta)\subset T_{\varepsilon_{0}}(\theta),$ the set $T_{\varepsilon_{0}%
}^{\prime}(\theta)$ is closed and so compact, by (\ref{di}) and the current
hypothesis (i). We apply Theorem \ref{thmcompact1} to\ get
\begin{equation}
\partial f(\theta)=\partial\ell(\theta)=\partial(\operatorname*{cl}%
\ell)(\theta)=\bigcap\limits_{L\in\mathcal{F}(\theta)}\overline
{\operatorname*{co}}\left\{  \bigcup\limits_{T^{\prime}(\theta)}%
\partial((\operatorname*{cl}\ell_{t})+\mathrm{I}_{\overline{L\cap
\operatorname*{dom}(\operatorname*{cl}\ell)}})(\theta)\right\}  ,\label{th}%
\end{equation}
where (recall (\ref{di}))
\begin{equation}
T^{\prime}(\theta)=\left\{  t\in T\mid(\operatorname*{cl}f_{t})(\theta
)=0\right\}  \subset\left\{  t\in T\mid f_{t}(\theta)=0\right\}  .\label{rer2}%
\end{equation}
Since, for every $L\in\mathcal{F}(\theta),$\
\[
(\operatorname*{cl}\ell_{t})+\mathrm{I}_{\overline{L\cap\operatorname*{dom}%
(\operatorname*{cl}\ell)}}=\max\left\{  (\operatorname*{cl}f_{t}%
)+\mathrm{I}_{\overline{L\cap\operatorname*{dom}(\operatorname*{cl}\ell)}%
},-2\varepsilon_{0}\right\}  ,
\]
for all $t\in T^{\prime}(\theta)$ we have
\[
\partial((\operatorname*{cl}\ell_{t})+\mathrm{I}_{\overline{L\cap
\operatorname*{dom}(\operatorname*{cl}\ell)}})(\theta)=\partial
((\operatorname*{cl}f_{t})+\mathrm{I}_{\overline{L\cap\operatorname*{dom}%
(\operatorname*{cl}\ell)}})(\theta)\subset\partial((\operatorname*{cl}%
f_{t})+\mathrm{I}_{\overline{L\cap\operatorname*{dom}f}})(\theta)\text{,}%
\]
and similarly,
\[
\partial((\operatorname*{cl}f_{t})+\mathrm{I}_{\overline{L\cap
\operatorname*{dom}f}})(\theta)\subset\partial(f_{t}+\mathrm{I}_{\overline
{L\cap\operatorname*{dom}f}})(\theta).
\]
Using successively the last two relations together with (\ref{rer2}), relation
(\ref{th}) implies
\begin{align}
\partial f(\theta) &  \subset\bigcap\limits_{L\in\mathcal{F}(\theta)}%
\overline{\operatorname*{co}}\left\{  \bigcup\limits_{t\in T^{\prime}(\theta
)}\partial((\operatorname*{cl}f_{t})+\mathrm{I}_{\overline{L\cap
\operatorname*{dom}f}})(\theta)\right\}  \label{toeps}\\
&  \subset\bigcap\limits_{L\in\mathcal{F}(\theta)}\overline{\operatorname*{co}%
}\left\{  \bigcup\limits_{t\in T^{\prime}(\theta)}\partial((\operatorname*{cl}%
f_{t})+\mathrm{I}_{L\cap\operatorname*{dom}f})(\theta)\right\}  \label{toeps2}%
\\
&  \subset\bigcap\limits_{L\in\mathcal{F}(\theta)}\overline{\operatorname*{co}%
}\left\{  \bigcup\limits_{t\in T(\theta)}\partial(f_{t}+\mathrm{I}%
_{L\cap\operatorname*{dom}f})(\theta)\right\}  ,\nonumber
\end{align}
which gives the desired\ inclusion"$\supset$".
\end{dem}

\begin{rem}
\label{rqq}\emph{As (\ref{toeps}) shows, we have proved the following
equivalent formula (under the assumptions of Theorem \ref{thmcompact})\ }%
\begin{align}
\partial f(x)  &  =\bigcap\limits_{L\in\mathcal{F}(x)}\overline
{\operatorname*{co}}\left\{  \bigcup\limits_{t\in T^{\prime}(x)}%
\partial((\operatorname*{cl}f_{t})+\mathrm{I}_{\overline{L\cap
\operatorname*{dom}f}})(x)\right\} \label{rqq1}\\
&  =\bigcap\limits_{L\in\mathcal{F}(x)}\overline{\operatorname*{co}}\left\{
\bigcup\limits_{t\in T^{\prime}(x)}\partial((\operatorname*{cl}f_{t}%
)+\mathrm{I}_{L\cap\operatorname*{dom}f})(x)\right\}  , \label{rqq2}%
\end{align}
\emph{where} $T^{\prime}(x)=\left\{  t\in T\mid(\operatorname*{cl}%
f_{t})(x)=f(x)\right\}  .$
\end{rem}

\bigskip

In order to characterize the subdifferential of the supremum function $f$ by
using only the functions $f_{t},$ rather than the augmented ones
$f_{t}+\mathrm{I}_{L\cap\operatorname*{dom}f},$ we provide another formula in
the following theorem.

\begin{theo}
\label{spe1}Let $\left\{  f_{t},t\in T\right\}  $ and $x\in X$ be as in
Theorem \emph{\ref{thmcompact}}. Then
\[
\partial f(x)=\bigcap_{\substack{\varepsilon>0\\L\in\mathcal{F}(x)}%
}\overline{\operatorname*{co}}\left\{  \bigcup\limits_{t\in T(x)}%
\partial_{\varepsilon}f_{t}(x)+\mathrm{N}_{L\cap\operatorname*{dom}%
f}(x)\right\}  .
\]

\end{theo}

\begin{dem}
The inclusion "$\supset$" is easy to prove. Indeed, we have that, for all
$t\in T(x)$ and $\varepsilon>0,$
\[
\partial_{\varepsilon}f_{t}(x)\subset\partial_{\varepsilon}f(x)\text{,}%
\]
and so, for all $L\in\mathcal{F}(x),$
\[
\partial_{\varepsilon}f_{t}(x)+\mathrm{N}_{L\cap\operatorname*{dom}%
f}(x)\subset\partial_{\varepsilon}f(x)+\mathrm{N}_{L\cap\operatorname*{dom}%
f}(x)\subset\partial_{\varepsilon}(f+\mathrm{I}_{L\cap\operatorname*{dom}%
f})(x)=\partial_{\varepsilon}(f+\mathrm{I}_{L})(x).
\]
Hence,
\begin{align*}
\bigcap_{\varepsilon>0,\text{ }L\in\mathcal{F}(x)}\overline{\operatorname*{co}%
}\left\{  \bigcup\limits_{t\in T(x)}\partial_{\varepsilon}f_{t}(x)+\mathrm{N}%
_{L\cap\operatorname*{dom}f}(x)\right\}   &  \subset\bigcap_{L\in
\mathcal{F}(x)}\bigcap_{\varepsilon>0}\partial_{\varepsilon}(f+\mathrm{I}%
_{L})(x)\\
&  =\bigcap_{L\in\mathcal{F}(x)}\partial(f+\mathrm{I}_{L})(x)=\partial f(x).
\end{align*}

Let us prove the inclusion "$\subset$". As in the previous theorems we suppose
that $x=\theta,$ $\partial f(\theta)\neq\emptyset,$ and $f(\theta
)=(\operatorname*{cl}f)(\theta)=0.$ Then, by (\ref{rqq1}),
\begin{equation}
\partial f(\theta)=\bigcap\limits_{L\in\mathcal{F}(\theta)}\overline
{\operatorname*{co}}\left\{  \bigcup\limits_{t\in T^{\prime}(\theta)}%
\partial((\operatorname*{cl}f_{t})+\mathrm{I}_{\overline{L\cap
\operatorname*{dom}f}})(\theta)\right\}  ,\label{rw}%
\end{equation}
where $T^{\prime}(\theta)=\left\{  t\in T\mid(\operatorname*{cl}f_{t}%
)(\theta)=0\right\}  .$ We are going to apply \cite[Theorem 12]{vv} and to
this purpose we see, from the one hand, that for every $t\in T$\
\[
\emptyset\neq\operatorname*{ri}(\overline{L\cap\operatorname*{dom}%
f})=\operatorname*{ri}(L\cap\operatorname*{dom}f)\subset\operatorname*{dom}%
f_{t}\subset\operatorname*{dom}(\operatorname*{cl}f_{t}),
\]
and this entails
\[
\operatorname*{dom}(\operatorname*{cl}f_{t})\cap\operatorname*{ri}%
(\operatorname*{dom}\mathrm{I}_{\overline{L\cap\operatorname*{dom}f}}%
)\neq\emptyset.
\]
On the other hand, the restriction of the function $\mathrm{I}_{\overline
{L\cap\operatorname*{dom}f}}$ to the affine hull of its domain is continuous
on $\operatorname*{ri}(\operatorname*{dom}\mathrm{I}_{\overline{L\cap
\operatorname*{dom}f}}),$ and so, \cite[Theorem 12]{vv} applies and yields
\begin{equation}
\partial((\operatorname*{cl}f_{t})+\mathrm{I}_{\overline{L\cap
\operatorname*{dom}f}})(\theta)=\bigcap\limits_{\varepsilon>0}%
\operatorname*{cl}\left(  \partial_{\varepsilon}(\operatorname*{cl}%
f_{t})(\theta)+\mathrm{N}_{\overline{L\cap\operatorname*{dom}f}}%
(\theta)\right)  .\label{sf}%
\end{equation}
If\ $t\in T^{\prime}(\theta)$ we have $(\operatorname*{cl}f_{t})(\theta
)=f_{t}(\theta)=0,$ and for all $\varepsilon>0$
\[
\partial_{\varepsilon}(\operatorname*{cl}f_{t})(\theta)+\mathrm{N}%
_{\overline{L\cap\operatorname*{dom}f}}(\theta)\subset\partial_{\varepsilon
}f_{t}(\theta)+\mathrm{N}_{\overline{L\cap\operatorname*{dom}f}}%
(\theta)=\partial_{\varepsilon}f_{t}(\theta)+\mathrm{N}_{L\cap
\operatorname*{dom}f}(\theta),
\]
so that (\ref{sf}) yields
\[
\partial((\operatorname*{cl}f_{t})+\mathrm{I}_{\overline{L\cap
\operatorname*{dom}f}})(\theta)\subset\operatorname*{cl}\left(  \partial
_{\varepsilon}f_{t}(\theta)+\mathrm{N}_{L\cap\operatorname*{dom}f}%
(\theta)\right)  .
\]
Consequently, by (\ref{rw}),
\begin{align*}
\partial f(\theta) &  \subset\bigcap\limits_{L\in\mathcal{F}(\theta)}%
\overline{\operatorname*{co}}\left\{  \bigcup\limits_{t\in T^{\prime}(\theta
)}\operatorname*{cl}\left(  \partial_{\varepsilon}f_{t}(\theta)+\mathrm{N}%
_{L\cap\operatorname*{dom}f}(\theta)\right)  \right\}  \\
&  =\bigcap\limits_{L\in\mathcal{F}(\theta)}\overline{\operatorname*{co}%
}\left\{  \bigcup\limits_{t\in T^{\prime}(\theta)}\partial_{\varepsilon}%
f_{t}(\theta)+\mathrm{N}_{L\cap\operatorname*{dom}f}(\theta)\right\}  .
\end{align*}
Since $T^{\prime}(\theta)\subset T(\theta),$ the aimed inclusion follows by
intersecting over $\varepsilon>0.$
\end{dem}

\bigskip

To avoid in Theorem \ref{thmcompact} the intersection\ over sets
$L\in\mathcal{F}(x),$ one has to require extra conditions relying either on
the space $X$ or on the function $f.$ We start with the following result,
whose first part\ is similar to the finite-dimensional-like\ result
established in Theorem \ref{thmcompact0}.

\begin{cor}
\label{sep2}Let $\left\{  f_{t},t\in T\right\}  $ and $x\in X$ be as in
Theorem \emph{\ref{thmcompact}}. Suppose that $\operatorname*{ri}%
(\operatorname*{dom}f)\neq\emptyset$ and $f_{\mid\operatorname*{aff}%
(\operatorname*{dom}f)}$ is continuous on $\operatorname*{ri}%
(\operatorname*{dom}f).$ Then
\[
\partial f(x)=\overline{\operatorname*{co}}\left\{  \bigcup\limits_{t\in
T(x)}\partial(f_{t}+\mathrm{I}_{\operatorname*{dom}f})(x)\right\}  .
\]
Moreover, if the following two conditions hold for all $t\in T(x):\smallskip$

\emph{(a)} $\operatorname*{ri}(\operatorname*{dom}f_{t})\cap
\operatorname*{dom}f\neq\emptyset,$

\emph{(b)} $(f_{t})_{\mid\operatorname*{aff}(\operatorname*{dom}f_{t})}$ is
continuous on $\operatorname*{ri}(\operatorname*{dom}f_{t}),$\newline then
\[
\partial f(x)=\overline{\operatorname*{co}}\left\{  \bigcup\limits_{t\in
T(x)}\partial f_{t}(x)+\mathrm{N}_{\operatorname*{dom}f}(x)\right\}  .
\]

\end{cor}

\begin{dem}
The inclusion "$\supset$" comes from Theorem \ref{thmcompact}, due to the
following relation which is true for every $L\in\mathcal{F}(x)$ and $t\in
T(x),$%
\[
\partial(f_{t}+\mathrm{I}_{\operatorname*{dom}f})(x)\subset\partial
(f_{t}+\mathrm{I}_{L\cap\operatorname*{dom}f})(x),
\]
so that
\[
\overline{\operatorname*{co}}\left\{  \bigcup\limits_{t\in T(x)}\partial
(f_{t}+\mathrm{I}_{\operatorname*{dom}f})(x)\right\}  \subset\bigcap
\limits_{L\in\mathcal{F}(x)}\overline{\operatorname*{co}}\left\{
\bigcup\limits_{t\in T(x)}\partial(f_{t}+\mathrm{I}_{L\cap\operatorname*{dom}%
f})(x)\right\}  =\partial f(x).
\]

To prove the inclusion "$\subset$" we may suppose that $x\in
\operatorname*{dom}(\partial f).$

First, we assume that the functions $f_{t}+\mathrm{I}_{L\cap
\operatorname*{dom}f},$ $t\in T$ are proper. Take $t\in T(x).$ By the current
assumption, we choose a point $x_{0}\in\operatorname*{ri}(\operatorname*{dom}%
f).$ Fix an $L\in\mathcal{F}(x)$ such that $x_{0}\in L.$ If
\[
q_{t}:=f_{t}+\mathrm{I}_{\operatorname*{dom}f},
\]
then $\operatorname*{dom}q_{t}=\operatorname*{dom}f$ and
\[
\operatorname*{ri}(\operatorname*{dom}q_{t})\cap\operatorname*{ri}%
(\operatorname*{dom}\mathrm{I}_{L})=\operatorname*{ri}(\operatorname*{dom}%
f)\cap L\neq\emptyset,
\]
as $x_{0}\in\operatorname*{ri}(\operatorname*{dom}f)\cap L.$ Moreover, since
$q_{t\mid\operatorname*{aff}(\operatorname*{dom}q_{t})}\equiv q_{t\mid
\operatorname*{aff}(\operatorname*{dom}f)}\leq f_{\mid\operatorname*{aff}%
(\operatorname*{dom}f)},$ the (proper convex) function $q_{t\mid
\operatorname*{aff}(\operatorname*{dom}q_{t})}$ is continuous on
$\operatorname*{ri}(\operatorname*{dom}q_{t}).$ It is clear that
$\mathrm{I}_{L\mid\operatorname*{aff}(\operatorname*{dom}\mathrm{I}_{L})}$
$(\equiv0_{\mid L})$ is continuous on $L.$ Then from \cite[Theorem 5]{SIOPT16}
it follows that
\begin{align}
\partial(f_{t}+\mathrm{I}_{L\cap\operatorname*{dom}f})(x)  &  =\partial
(q_{t}+\mathrm{I}_{L})(x)\nonumber\\
&  =\operatorname*{cl}(\partial q_{t}(x)+\partial\mathrm{I}_{L}%
(x))=\operatorname*{cl}(\partial(f_{t}+\mathrm{I}_{\operatorname*{dom}%
f})(x)+L^{\bot}). \label{ta}%
\end{align}
Therefore, by Theorem \ref{thmcompact} we get
\begin{align*}
\partial f(x)  &  =\bigcap\limits_{L\in\mathcal{F}(x)}\overline
{\operatorname*{co}}\left\{  \bigcup\limits_{t\in T(x)}\partial(f_{t}%
+\mathrm{I}_{L\cap\operatorname*{dom}f})(x)\right\} \\
&  \subset\bigcap\limits_{L\in\mathcal{F}(x),\text{ }x_{0}\in L}%
\overline{\operatorname*{co}}\left\{  \bigcup\limits_{t\in T(x)}\partial
(f_{t}+\mathrm{I}_{L\cap\operatorname*{dom}f})(x)\right\} \\
&  =\bigcap\limits_{L\in\mathcal{F}(x),\text{ }x_{0}\in L}\overline
{\operatorname*{co}}\left\{  \bigcup\limits_{t\in T(x)}\operatorname*{cl}%
(\partial(f_{t}+\mathrm{I}_{\operatorname*{dom}f})(x)+L^{\bot})\right\}
\text{ \ (by (\ref{ta}))}\\
&  =\bigcap\limits_{L\in\mathcal{F}(x),\text{ }x_{0}\in L}\overline
{\operatorname*{co}}\left\{  \bigcup\limits_{t\in T(x)}\partial(f_{t}%
+\mathrm{I}_{\operatorname*{dom}f})(x)+L^{\bot}\right\}  \text{.}%
\end{align*}
Next, given a $U\in\mathcal{N}_{X^{\ast}}$ we choose an $F\in\mathcal{F}(x)$
such that $x_{0}\in F$ and $F^{\bot}\subset U.$ Hence, from the last inclusion
we obtain
\begin{align}
\partial f(x)  &  \subset\overline{\operatorname*{co}}\left\{  \bigcup
\limits_{t\in T(x)}\partial(f_{t}+\mathrm{I}_{\operatorname*{dom}%
f})(x)+F^{\bot}\right\} \nonumber\\
&  \subset\operatorname*{co}\left\{  \bigcup\limits_{t\in T(x)}\partial
(f_{t}+\mathrm{I}_{\operatorname*{dom}f})(x)\right\}  +F^{\bot}+U\nonumber\\
&  \subset\operatorname*{co}\left\{  \bigcup\limits_{t\in T(x)}\partial
(f_{t}+\mathrm{I}_{\operatorname*{dom}f})(x)\right\}  +2U, \label{argu}%
\end{align}
and the aimed inclusion follows by intersecting over $U\in\mathcal{N}%
_{X^{\ast}}.$

Now, to deal with the case when not necessarily all the $f_{t}$'s, $t\in T,$
are proper, we consider the functions%
\[
\hat{f}_{t}:=\max\{f_{t},f(x)-2\varepsilon_{0}\},\text{ }t\in T,\text{ }%
\hat{f}:=\max\hat{f}_{t}=\max\{f,f(x)-2\varepsilon_{0}\}.
\]
By Lemma \ref{lt}, and taking into account Remark \ref{remcl}, it follows that
the proper functions $\hat{f}_{t},$ $t\in T,$ satisfy the assumtpions of the
current theorem. So, from the previous part of the proof applied to the
$\hat{f}_{t}$'s we get
\begin{align}
\partial f(x)=\partial\hat{f}(x) &  =\overline{\operatorname*{co}}\left\{
\bigcup\limits_{\{t\in T\mid\hat{f}_{t}(x)=f(x)\}}\partial(\hat{f}%
_{t}+\mathrm{I}_{\operatorname*{dom}\hat{f}})(x)\right\}  \nonumber\\
&  =\overline{\operatorname*{co}}\left\{  \bigcup\limits_{t\in T(x)}%
\partial(\hat{f}_{t}+\mathrm{I}_{\operatorname*{dom}f})(x)\right\}
.\label{fin}%
\end{align}
Observe that for $t\in T(x)$ such that $\partial(\hat{f}_{t}+\mathrm{I}%
_{\operatorname*{dom}f})(x)\neq\emptyset$ it holds
\begin{align*}
f(x) &  =(\hat{f}_{t}+\mathrm{I}_{\operatorname*{dom}f})(x)=\operatorname*{cl}%
(\hat{f}_{t}+\mathrm{I}_{\operatorname*{dom}f})(x)\\
&  =\max\{\operatorname*{cl}(f_{t}+\mathrm{I}_{\operatorname*{dom}%
f})(x),f(x)-2\varepsilon_{0}\}\\
&  =\operatorname*{cl}(f_{t}+\mathrm{I}_{\operatorname*{dom}f})(x)\leq
(f_{t}+\mathrm{I}_{\operatorname*{dom}f})(x)\leq f(x),
\end{align*}
which implies that $\operatorname*{cl}(f_{t}+\mathrm{I}_{\operatorname*{dom}%
f})(x)=(f_{t}+\mathrm{I}_{\operatorname*{dom}f})(x)$ and so, the function
$f_{t}+\mathrm{I}_{\operatorname*{dom}f}$ is lsc at $x.$ Hence, $\partial
(\hat{f}_{t}+\mathrm{I}_{\operatorname*{dom}f})(x)=\partial(\max
\{f_{t}+\mathrm{I}_{\operatorname*{dom}f},f(x)-2\varepsilon_{0}\})(x)=\partial
(f_{t}+\mathrm{I}_{\operatorname*{dom}f})(x),$ and (\ref{fin}) implies
\[
\partial f(x)=\overline{\operatorname*{co}}\left\{  \bigcup\limits_{t\in
T(x)}\partial(f_{t}+\mathrm{I}_{\operatorname*{dom}f})(x)\right\}  ,
\]
and the proof the first assertion is finished.

The last assertion also comes from \cite[Theorem 5]{SIOPT16}. By the
accessibility lemma, condition (a) and the properties of $f$ imply that%
\[
\operatorname*{ri}(\operatorname*{dom}f_{t})\cap\operatorname*{ri}%
(\operatorname*{dom}f)\neq\emptyset\text{ \ for all }t\in T(x),
\]
and then
\begin{align*}
\partial f(x) &  =\overline{\operatorname*{co}}\left\{  \bigcup\limits_{t\in
T(x)}\partial(f_{t}+\mathrm{I}_{\operatorname*{dom}f})(x)\right\}  \\
&  =\overline{\operatorname*{co}}\left\{  \bigcup\limits_{t\in T(x)}%
\operatorname*{cl}(\partial f_{t}(x)+\mathrm{N}_{\operatorname*{dom}%
f}(x))\right\}  =\overline{\operatorname*{co}}\left\{  \bigcup\limits_{t\in
T(x)}\partial f_{t}(x)+\mathrm{N}_{\operatorname*{dom}f}(x)\right\}  .
\end{align*}

\end{dem}

By dropping conditions (a) and (b) in Corollary \ref{sep2} we derive the
following result.

\begin{cor}
\label{sep2b}Let $\left\{  f_{t},t\in T\right\}  $ and $x\in X$ be as in
Theorem \emph{\ref{thmcompact}}. Suppose that $\operatorname*{ri}%
(\operatorname*{dom}f)\neq\emptyset.$ Then
\[
\partial f(x)=\bigcap_{\varepsilon>0}\overline{\operatorname*{co}}\left\{
\bigcup\limits_{t\in T(x)}\partial_{\varepsilon}f_{t}(x)+\mathrm{N}%
_{\operatorname*{dom}f}(x)\right\}  .
\]

\end{cor}

\begin{dem}
The inclusion "$\supset$"\ follows from Theorem \ref{spe1} in view of the
following relation, for $t\in T(x),$ $L\in\mathcal{F}(x)$, and $\varepsilon
>0,$
\[
\partial_{\varepsilon}f_{t}(x)+\mathrm{N}_{\operatorname*{dom}f}%
(x)\subset\partial_{\varepsilon}f_{t}(x)+\mathrm{N}_{L\cap\operatorname*{dom}%
f}(x).
\]

To prove the inclusion "$\subset$" we suppose that $x\in\operatorname*{dom}%
(\partial f).$ By Theorem \ref{spe1} we have that
\begin{equation}
\partial f(x)=\bigcap\limits_{\varepsilon>0,\text{ }L\in\mathcal{F}%
(x)}\overline{\operatorname*{co}}\left\{  \bigcup\limits_{t\in T(x)}%
\partial_{\varepsilon}f_{t}(x)+\mathrm{N}_{L\cap\operatorname*{dom}%
f}(x)\right\}  .\label{kj}%
\end{equation}
The current assumption ensures, thanks to \cite[Theorem 5]{SIOPT16}, that for
every $L\in\mathcal{F}(x)$ such that $L\cap\operatorname*{ri}%
(\operatorname*{dom}f)\neq\emptyset$
\[
\mathrm{N}_{L\cap\operatorname*{dom}f}(x)=\operatorname*{cl}(\mathrm{N}%
_{\operatorname*{dom}f}(x)+L^{\bot}).
\]
Hence, by arguing as in (\ref{argu}), relation (\ref{kj}) leads to%
\begin{align*}
\partial f(x) &  \subset\bigcap\limits_{\substack{\varepsilon>0\\L\in
\mathcal{F}(x),\text{ }L\cap\operatorname*{ri}(\operatorname*{dom}%
f)\neq\emptyset}}\overline{\operatorname*{co}}\left\{  \bigcup\limits_{t\in
T(x)}\partial_{\varepsilon}f_{t}(x)+\operatorname*{cl}(\mathrm{N}%
_{\operatorname*{dom}f}(x)+L^{\bot})\right\}  \\
&  \subset\bigcap\limits_{\varepsilon>0,\text{ }U\in\mathcal{N}_{X^{\ast}}%
}\overline{\operatorname*{co}}\left\{  \bigcup\limits_{t\in T(x)}%
\partial_{\varepsilon}f_{t}(x)+\mathrm{N}_{\operatorname*{dom}f}(x)+U\right\}
\\
&  =\bigcap\limits_{\varepsilon>0}\overline{\operatorname*{co}}\left\{
\bigcup\limits_{t\in T(x)}\partial_{\varepsilon}f_{t}(x)+\mathrm{N}%
_{\operatorname*{dom}f}(x)\right\}  .
\end{align*}

\end{dem}

The following lemma is used in Theorem \ref{corcompcont}.

\begin{lem}
\label{lemf}Assume that $f=\sup_{t\in T}f_{t}$\ is finite and continuous at
some point. If\ $x\in X$ satisfies conditions $(i)$ and $(ii)$ of Theorem
\emph{\ref{thmcompact}}, then $\operatorname*{cl}f=\sup_{t\in T}%
(\operatorname*{cl}f_{t})$ and
\[
\partial f(x)=\overline{\operatorname*{co}}\left\{  \bigcup\limits_{t\in
T^{\prime}(x)}\partial(\operatorname*{cl}f_{t})(x)+\mathrm{N}%
_{\operatorname*{dom}f}(x)\right\}  ,
\]
where $T^{\prime}(x)=\left\{  t\in T\mid(\operatorname*{cl}f_{t}%
)(x)=f(x)\right\}  .$
\end{lem}

\begin{dem}
By Moreau-Rockafellar's Theorem (see, i.e., \cite[Proposition 10.3]%
{Moreaubook}), the continuity assumption ensures that, for every $t\in
T^{\prime}(x)$ such that $\partial((\operatorname*{cl}f_{t})+\mathrm{I}%
_{L\cap\operatorname*{dom}f})(x)\neq\emptyset,$ and every $L\in\mathcal{F}%
(x)\ $such that $L\cap\operatorname*{int}(\operatorname*{dom}f)\neq\emptyset,$%
\[
\partial((\operatorname*{cl}f_{t})+\mathrm{I}_{L\cap\operatorname*{dom}%
f})(x)=\partial(\operatorname*{cl}f_{t})(x)+\mathrm{N}_{\operatorname*{dom}%
f}(x)+L^{\bot},
\]
taking into account that the continuity of $f$ is inherited by the function
$\operatorname*{cl}f_{t},$ whose properness is a consequence of $\partial
((\operatorname*{cl}f_{t})+\mathrm{I}_{L\cap\operatorname*{dom}f}%
)(x)\neq\emptyset$. Observe that the equality $\operatorname*{cl}f=\sup_{t\in
T}(\operatorname*{cl}f_{t})$ follows according to Remark \ref{remcl}. Then,
using (\ref{rqq2}),\emph{ }%
\begin{align*}
\partial f(x) &  =\bigcap\limits_{L\in\mathcal{F}(x)}\overline
{\operatorname*{co}}\left\{  \bigcup\limits_{t\in T^{\prime}(x)}%
\partial((\operatorname*{cl}f_{t})+\mathrm{I}_{L\cap\operatorname*{dom}%
f})(x)\right\}  \\
&  \subset\bigcap\limits_{L\in\mathcal{F}(x),\text{ }L\cap\operatorname*{int}%
(\operatorname*{dom}f)\neq\emptyset}\overline{\operatorname*{co}}\left\{
\bigcup\limits_{t\in T^{\prime}(x)}\partial(\operatorname*{cl}f_{t}%
)(x)+\mathrm{N}_{\operatorname*{dom}f}(x)+L^{\bot}\right\}  \\
&  \subset\bigcap\limits_{U\in\mathcal{N}_{X^{\ast}}}\overline
{\operatorname*{co}}\left\{  \bigcup\limits_{t\in T^{\prime}(x)}%
\partial(\operatorname*{cl}f_{t})(x)+\mathrm{N}_{\operatorname*{dom}%
f}(x)+U\right\}  \\
&  =\overline{\operatorname*{co}}\left\{  \bigcup\limits_{t\in T^{\prime}%
(x)}\partial(\operatorname*{cl}f_{t})(x)+\mathrm{N}_{\operatorname*{dom}%
f}(x)\right\}  .
\end{align*}
Thus, as we can easily check that
\[
\overline{\operatorname*{co}}\left\{  \bigcup\limits_{t\in T^{\prime}%
(x)}\partial(\operatorname*{cl}f_{t})(x)+\mathrm{N}_{\operatorname*{dom}%
f}(x)\right\}  \subset\partial f(x),
\]
the conclusion follows.
\end{dem}

The following result simplifies Theorem \ref{thmcompact} when the supremum
function is continuous on the interior of its domain.

\begin{theo}
\label{corcompcont}Assume that the family of convex functions $\left\{
f_{t},t\in T\right\}  $ is such that $f$\ is finite and continuous at some
point. Let $x\in X$ satisfy conditions $(i)$ and $(ii)$ of Theorem
\emph{\ref{thmcompact}}. Then
\begin{align*}
\partial f(x) &  =\mathrm{N}_{\operatorname*{dom}f}(x)+\overline
{\operatorname*{co}}\left\{  \bigcup\limits_{t\in T(x)}\partial f_{t}%
(x)\right\}  \\
&  =\mathrm{N}_{\operatorname*{dom}f}(x)+\operatorname*{co}\left\{
\bigcup\limits_{t\in T(x)}\partial f_{t}(x)\right\}  \text{ \ }(\text{when
}X=\mathbb{R}^{n}).
\end{align*}

\end{theo}

\begin{dem}
The inclusion "$\supset$" is straighforward, and thus we only need to check
the converse inclusion in the nontrivial case $\partial f(x)\neq\emptyset;$
therefore, we may assume that $x=\theta$ and $f(\theta)=0.$ From Lemma
\ref{lemf} we have that
\[
\partial f(\theta)=\overline{\operatorname*{co}}\left\{  \bigcup\limits_{t\in
T^{\prime}(\theta)}\partial(\operatorname*{cl}f_{t})(\theta)+\mathrm{N}%
_{\operatorname*{dom}f}(\theta)\right\}  ,
\]
where $T^{\prime}(x)=\left\{  t\in T\mid(\operatorname*{cl}f_{t}%
)(x)=f(x)\right\}  .$ By the current continuity assumption we choose $x_{0}\in
X,$ $W\in\mathcal{N}_{X}$ and $m\geq0$ such that for all $w\in W$
\begin{equation}
f(x_{0}+w)\leq m;\label{cu}%
\end{equation}
hence,
\begin{align*}
\mathrm{\sigma}_{\cup_{t\in T^{\prime}(\theta)}\partial(\operatorname*{cl}%
f_{t})(\theta)}(x_{0}+w) &  \leq\sup_{t\in T^{\prime}(\theta)}%
((\operatorname*{cl}f_{t})(x_{0}+w)-(\operatorname*{cl}f_{t})(\theta))\\
&  \leq\sup_{t\in T^{\prime}(\theta)}f_{t}(x_{0}+w)\leq f(x_{0}+w)\leq m,
\end{align*}
showing that the proper lsc convex function $\mathrm{\sigma}_{\cup_{t\in
T^{\prime}(\theta)}\partial(\operatorname*{cl}f_{t})(\theta)}$ is also
continuous at $x_{0}.$ Thus, by applying once again Moreau-Rockafellar's
Theorem, since \textrm{$\sigma$}$_{\mathrm{N}_{\operatorname*{dom}f}(\theta
)}(x_{0})\leq0$ we obtain
\begin{align*}
\overline{\operatorname*{co}}\left\{  \bigcup\limits_{t\in T^{\prime}(\theta
)}\partial(\operatorname*{cl}f_{t})(\theta)+\mathrm{N}_{\operatorname*{dom}%
f}(\theta)\right\}   &  =\partial\left(  \mathrm{\sigma}_{\cup_{t\in
T^{\prime}(\theta)}\partial(\operatorname*{cl}f_{t})(\theta)+\mathrm{N}%
_{\operatorname*{dom}f}(\theta)}\right)  (\theta)\\
&  =\partial\left(  \mathrm{\sigma}_{\cup_{t\in T^{\prime}(\theta)}%
\partial(\operatorname*{cl}f_{t})(\theta)}+\mathrm{\sigma}_{\mathrm{N}%
_{\operatorname*{dom}f}(\theta)}\right)  (\theta)\\
&  =\partial\mathrm{\sigma}_{\cup_{t\in T^{\prime}(\theta)}\partial
(\operatorname*{cl}f_{t})(\theta)}(\theta)+\partial\mathrm{\sigma}%
_{\mathrm{N}_{\operatorname*{dom}f}(\theta)}(\theta).
\end{align*}
Using again the well-known formula of the subdifferential of the support
function, we obtain
\begin{align}
\overline{\operatorname*{co}}\left\{  \bigcup\limits_{t\in T^{\prime}(\theta
)}\partial(\operatorname*{cl}f_{t})(\theta)+\mathrm{N}_{\operatorname*{dom}%
f}(\theta)\right\}    & =\overline{\operatorname*{co}}\left\{  \bigcup
\limits_{t\in T^{\prime}(\theta)}\partial(\operatorname*{cl}f_{t}%
)(\theta)\right\}  +\mathrm{N}_{\operatorname*{dom}f}(\theta)\label{final}\\
& \subset\overline{\operatorname*{co}}\left\{  \bigcup\limits_{t\in T(\theta
)}\partial f_{t}(\theta)\right\}  +\mathrm{N}_{\operatorname*{dom}f}%
(\theta),\nonumber
\end{align}
and the first formula follows.

To prove the second statment of the theorem, when $X=\mathbb{R}^{n},$ observe
by (\ref{final}) that
\[
\partial f(\theta)=\overline{\operatorname*{co}}\left\{  \bigcup\limits_{t\in
T^{\prime}(\theta)}\partial(\operatorname*{cl}f_{t})(\theta)\right\}
+\mathrm{N}_{\operatorname*{dom}f}(\theta),
\]
and thus, we shall show that
\[
\overline{\operatorname*{co}}\left\{  \bigcup\limits_{t\in T^{\prime}(\theta
)}\partial(\operatorname*{cl}f_{t})(\theta)\right\}  \subset\operatorname*{co}%
\left\{  \bigcup\limits_{t\in T^{\prime}(\theta)}\partial(\operatorname*{cl}%
f_{t})(\theta)\right\}  +\mathrm{N}_{\operatorname*{dom}f}(\theta).
\]
Fix $u$ in the left-hand side. By taking into account Caratheodory's Theorem
we choose
\[
(\lambda_{1,i},\cdots,\lambda_{n+1,i})\in\Delta_{n+1}:=\left\{  (\lambda
_{1},\cdots,\lambda_{n+1})\mid\lambda_{i}\geq0,\sum\nolimits_{i=1}%
^{n+1}\lambda_{i}=1\right\}  ,
\]
$t_{1,i},\cdots,t_{n+1,i}\in T^{\prime}(\theta),$ and $y_{1,i}\in
\partial(\operatorname*{cl}f_{t_{1,i}})(\theta),$ $\cdots,y_{n+1,i}\in
\partial(\operatorname*{cl}f_{t_{n+1,i}})(\theta),$ $i=1,2,\cdots,$ such that
\begin{equation}
y_{i}:=\lambda_{1,i}y_{1,i}+\cdots+\lambda_{n+1,i}y_{n+1,i}\rightarrow
u.\label{uu}%
\end{equation}
We may assume that for all $k=1,\cdots,n+1,$ there are subnets $(\lambda
_{k,\alpha})_{\alpha\in\Upsilon}$ and $(t_{k,\alpha})_{\alpha\in\Upsilon}$ of
$(\lambda_{k,i})_{i}$ and $(t_{k,i})_{i}$, respectively, such that
$\lambda_{k,\alpha}\rightarrow\lambda_{k}\in\left[  0,1\right]  $ and
$t_{k,\alpha}\rightarrow t_{k}\in T^{\prime}(\theta)$ (remember that
$T^{\prime}(\theta)$ is a closed subset of the compact set $T_{\varepsilon
_{0}}(\theta)$)$.$ Fix $\alpha\in\Upsilon$ and $k\in\left\{  1,\cdots
,n+1\right\}  .$ Then, using (\ref{cu}), for\ every $w\in W$%
\begin{align}
\left\langle \lambda_{k,\alpha}y_{k,\alpha},x_{0}+w\right\rangle  &
\leq\lambda_{k,\alpha}((\operatorname*{cl}f_{t_{k,\alpha}})(x_{0}%
+w)-(\operatorname*{cl}f_{t_{k,\alpha}})(\theta))\nonumber\\
& \leq\lambda_{k,\alpha}f(x_{0}+w)\leq\lambda_{k,\alpha}m\leq m,\label{la}%
\end{align}
and, taking $w=\theta,$ $\left\langle \lambda_{k,\alpha}y_{k,\alpha}%
,x_{0}\right\rangle \leq m.$ Due to (\ref{uu}), implying that
\[
\left\langle \lambda_{1,\alpha}y_{1,\alpha}+\cdots+\lambda_{n+1,\alpha
}y_{n+1,\alpha},x_{0}\right\rangle \rightarrow\left\langle u,x_{0}%
\right\rangle ,
\]
we deduce that the net $\left(  \left\langle \lambda_{k,\alpha}y_{k,\alpha
},x_{0}\right\rangle \right)  _{\alpha\in\Upsilon}$ is bounded. Consequently,
(\ref{la}) ensures that for some $r>0$
\[
\left\langle \lambda_{k,\alpha}y_{k,\alpha},w\right\rangle \leq r\text{ \ for
all }w\in W,
\]
showing that $(\lambda_{k,\alpha}y_{k,\alpha})_{\alpha\in\Upsilon}\subset
rW^{\circ}.$ Without loss of generality, we may suppose that $(\lambda
_{k,\alpha}y_{k,\alpha})_{\alpha\in\Upsilon}$ converges to some $y_{k}%
\in\mathbb{R}^{n}.$

Let%
\[
K_{+}:=\{k=1,\cdots,n+1\mid\lambda_{k}>0\}\text{ and }K_{0}:=\{k=1,\cdots
,n+1\mid\lambda_{k}=0\}.
\]
If\ $k\in K_{+},$ from\ $y_{k,\alpha}\in\partial(\operatorname*{cl}%
f_{t_{k,\alpha}})(\theta)$ we obtain, for all $z\in\operatorname*{dom}f,$
\[
\left\langle y_{k,\alpha},z\right\rangle \leq(\operatorname*{cl}%
f_{t_{k,\alpha}})(z)-(\operatorname*{cl}f_{t_{k,\alpha}})(\theta
)=(\operatorname*{cl}f_{t_{k,\alpha}})(z),\text{ }%
\]
which, by assumption (ii) and after passing to the limit on $\alpha\in
\Upsilon,$ yields
\[
\left\langle \lambda_{k}^{-1}y_{k},z\right\rangle \leq\limsup_{\alpha
\in\Upsilon}(\operatorname*{cl}f_{t_{k,\alpha}})(z)\leq(\operatorname*{cl}%
f_{t_{k}})(z)\leq f_{t_{k}}(z)\text{ \ for all }z\in\operatorname*{dom}%
(\operatorname*{cl}f)\text{ }(\supset\operatorname*{dom}f).
\]
By taking into account Moreau-Rockafellar's Theorem,\ this shows that
\begin{equation}
\lambda_{k}^{-1}y_{k}\in\partial(f_{t_{k}}+\mathrm{I}_{\operatorname*{dom}%
f})(\theta)=\partial f_{t_{k}}(\theta)+\mathrm{N}_{\operatorname*{dom}%
f}(\theta).\label{hh}%
\end{equation}
If $k\in K_{0},$ for all $z\in\operatorname*{dom}f$ we have
\[
\left\langle \lambda_{k,\alpha}y_{k,\alpha},z\right\rangle \leq\lambda
_{k,\alpha}((\operatorname*{cl}f_{t_{k,\alpha}})(z)-(\operatorname*{cl}%
f_{t_{k,\alpha}})(\theta))\leq\lambda_{k,\alpha}f(z),
\]
and by taking the limit on $\alpha\in\Upsilon$ we obtain that $y_{k}%
\in\mathrm{N}_{\operatorname*{dom}f}(x).$ This, together with (\ref{hh}) and
(\ref{uu}), leads us to
\[
u=\sum\limits_{k\in K_{+}}\lambda_{k}(\lambda_{k}^{-1}y_{k})+\sum\limits_{k\in
K_{0}}y_{k}\in\operatorname*{co}\left\{  \bigcup\limits_{t\in T(x)}\partial
f_{t}(\theta)\right\}  +\mathrm{N}_{\operatorname*{dom}f}(\theta),
\]
as the inclusion "$\subset$" follows.
\end{dem}

Finally, we give Valadier's formula under slightly weaker conditions.

\begin{cor}
Assume that the family of convex functions $\left\{  f_{t},t\in T\right\}  $
is such that $f$\ is finite and continuous at $x\in X.$ Suppose that for some
$\varepsilon_{0}>0$\emph{:}

$(i)$ the set $T_{\varepsilon_{0}}(x)$ is compact,

$(ii)$ the functions $t\mapsto f_{t}(z),$ $z\in\operatorname*{dom}f,$\ are usc
on $T_{\varepsilon_{0}}(x).$ Then
\begin{align*}
\partial f(x) &  =\overline{\operatorname*{co}}\left\{  \bigcup\limits_{t\in
T(x)}\partial f_{t}(x)\right\}  \\
&  =\operatorname*{co}\left\{  \bigcup\limits_{t\in T(x)}\partial
f_{t}(x)\right\}  \text{ \ }(\text{when }X=\mathbb{R}^{n}).
\end{align*}

\end{cor}

\begin{dem}
Since the inclusion "$\supset$" is immediate, we only need to show
the\ converse one when $\partial f(x)\neq\emptyset;$ hence, we may suppose
that $x=\theta$ and $f(\theta)=0.$ We choose an open $\theta$-neighborhood
$U\subset X$ and an $m\geq0$ such that
\[
-2\varepsilon_{0}<f(u)\leq m\text{ \ for all }u\in U.
\]
We denote,
\[
\tilde{f}_{t}:=\max\{f_{t},-2\varepsilon_{0}\}+\mathrm{I}_{U},\text{ }t\in
T,\text{ and }\tilde{f}:=\sup_{t\in T}\tilde{f}_{t}=\max\{f,-2\varepsilon
_{0}\}+\mathrm{I}_{U}.
\]
It is clear that the proper convex function\ $\tilde{f}$ is finite and
continuous at $\theta,$ with $\tilde{f}(\theta)=0,$ and that
\begin{equation}
\{t\in T\mid\tilde{f}_{t}(\theta)\geq\tilde{f}(\theta)-\varepsilon
\}=T_{\varepsilon}(\theta)\text{ \ for all }\varepsilon\in\left[
0,\varepsilon_{0}\right]  .\label{inf}%
\end{equation}
Also, from the inequalities
\[
\tilde{f}_{t}(u)\leq\tilde{f}(u)\leq\max\{m,-2\varepsilon_{0}\}=m\ \ \text{for
all }t\in T\text{ and }u\in U,
\]
we deduce that all the proper convex functions $\tilde{f}_{t}$, $t\in T,$ and
$\tilde{f}$ are continuous on $U.$ So,
\[
\tilde{f}_{t}(u)=(\operatorname*{cl}\tilde{f}_{t})(u)\text{, }t\in T,\text{
and }\tilde{f}(u)=(\operatorname*{cl}\tilde{f})(u)\text{ \ for all }u\in U.
\]
Therefore, by assumption (ii), for every $z\in\operatorname*{dom}%
(\operatorname*{cl}\tilde{f})$ $(=\operatorname*{dom}(\operatorname*{cl}f)\cap
U\subset U\subset\operatorname*{dom}f)$ the function
\[
t\mapsto(\operatorname*{cl}\tilde{f}_{t})(z)=\tilde{f}_{t}(z)=\max
\{f_{t}(z),-2\varepsilon_{0}\}
\]
is usc on $T_{\varepsilon_{0}}(\theta).$ Consequently, by Theorem
\ref{corcompcont} applied to the family $\{\tilde{f}_{t},$ $t\in T\}$ we
obtain (recall (\ref{inf}))
\begin{equation}
\partial f(\theta)=\partial\tilde{f}(\theta)=\overline{\operatorname*{co}%
}\left\{  \bigcup\limits_{t\in T(\theta)}\partial\tilde{f}_{t}(\theta
)\right\}  .\label{yut}%
\end{equation}
Moreover, if $t\in T(\theta)$ is such that $\partial\tilde{f}_{t}(\theta
)\neq\emptyset,$ then $\tilde{f}_{t}$ is lsc at $\theta$, and so using Remark
\ref{remcl} (and (\ref{inf})),
\begin{align*}
0 &  =\tilde{f}(\theta)=\tilde{f}_{t}(\theta)=(\operatorname*{cl}\tilde{f}%
_{t})(\theta)=(\operatorname*{cl}(\max\{f_{t},-2\varepsilon_{0}\}))(\theta)\\
&  =\max\{(\operatorname*{cl}f_{t})(\theta),-2\varepsilon_{0}\}\leq\max
\{f_{t}(\theta),-2\varepsilon_{0}\}\leq\max\{f(\theta),-2\varepsilon_{0}\}=0.
\end{align*}
Hence, $(\operatorname*{cl}f_{t})(\theta)=f_{t}(\theta)=0\ $and so $f_{t}$ is
lsc at $\theta.$ This implies that $\partial\tilde{f}_{t}(\theta)=\partial
f_{t}(\theta),$ and relation (\ref{yut}) gives
\[
\partial f(\theta)=\overline{\operatorname*{co}}\left\{  \bigcup\limits_{t\in
T(x)}\partial f_{t}(x)\right\}  .
\]

The second assertion is obtained\ in the same way.
\end{dem}


\begin{thebibliography}{99}                                                                                               %


\bibitem {Borwlewis}\textsc{Borwein, J. M. and Lewis,} \textsc{A.S.\ }%
\emph{Partially finite convex programming. I. Quasi relative interiors and
duality theory}. Math. Programming 57 (1992), Ser. B, 15--48.

\bibitem {Brondsted72}\textsc{Brøndsted}, \textsc{A. }\emph{On the
subdifferential of the supremum of two convex functions}, Math. Scand. 31
(1972), 225--230.

\bibitem {RcAhMl16}\textsc{Correa, R., Hantoute, A. and López, M. A.
}\emph{Valadier-like formulas for the supremum function I, }submitted, 2016.

\bibitem {SIOPT16}\textsc{Correa, R., Hantoute, A. and López, M. A.
}\emph{Towards Supremum-Sum subdifferential calculus free of qualification
conditions. }SIAM J. Optim. 26 (2016), 2219--2234.

\bibitem {vv}\textsc{Correa, R., Hantoute, A. and López, M. A. }\emph{Weaker
conditions for subdifferential calculus of convex functions.} J. Funct. Anal.
271 (2016), 1177--1212.

\bibitem {Han06}\textsc{Hantoute, A.} \emph{Subdifferential set of the
supremum of lower semi-continuous convex functions and the conical hull
intersection property}. Top 14 (2006), 355--374.

\bibitem {HanLop08}\textsc{Hantoute, A. and López, M. A.} \emph{A complete
characterization of the subdifferential set of the supremum of an arbitrary
family of convex functions}. J. Convex Anal. 15 (2008), 831--858.

\bibitem {WELL}\textsc{Hantoute, A. and \textsc{L}ópez, M. A.}
\emph{Characterization of total ill-posedness in linear semi-infinite
optimization}. J. Comput. Appl. Math. 217 (2008), 350--364.

\bibitem {HanLopZal2008}\textsc{Hantoute, A., López and M. A., Z\u{a}linescu,
C.} \emph{Subdifferential calculus rules in convex analysis: a unifying
approach via pointwise supremum functions.} SIAM J. Optim. 19 (2008), 863--882.

\bibitem {Ioffe11}\textsc{Ioffe, A. D.}, \emph{A note on subdifferentials of
pointwise suprema. }Top 20 (2012), 456-466.

\bibitem {Ioffe72}\textsc{Ioffe,A.~D. and Levin, V.\ L. }%
\emph{Subdifferentials of convex functions}, Trudy Mos. Mat. Obs., 26 (1972),
3--73 (Russian).

\bibitem {LiNg11}\textsc{Li, C. and Ng, K. F. }\emph{Subdifferential Calculus
Rules for Supremum Functions in Convex Analysis.} SIAM J. Optim. 21 (2011), 782--797.

\bibitem {Olopez}\textsc{López, O. and Thibault, L, }\emph{Sequential formula
for subdifferential of upper envelope of convex functions}. J. Nonlinear
Convex Anal. 14 (2013), 377--388.

\bibitem {MOR}\textsc{Mordukhovich, B and Nghia, T. T. A.}
\emph{Subdifferentials of nonconvex supremum functions and their applications
to semi-infinite and infinite programs with Lipschitzian data}. SIAM J. Optim.
23 (2013), 406--431. 

\bibitem {Moreaubook}\textsc{Moreau, J. J.}\ Fonctionnelles convexes, in:
Séminaire sur les équations aux dérivées partielles. Collège de France, 1966.

\bibitem {Rockafellar66}\textsc{Rockafellar, R. T. }\emph{Characterization of
the subdifferentials of convex functions}. Pacific J. Math. 17 (1966), 497--510.

\bibitem {Soloviev01}\textsc{Solov'ev, V.~N.}, \emph{The subdifferential and
the directional derivatives of the maximum of a family of convex functions},
Izv. Ross Akad. Nauk Ser. Mat., 65 (2001), 107--132.

\bibitem {Valadier}\textsc{Valadier, V. }\emph{Sous-différentiels d'une borne
supérieure et d'une somme continue de fonctions convexes}, C. R. Acad. Sci.
Paris Sér. A-B Math., 268 (1969), 39--42.

\bibitem {ZalinescuBook}\textsc{Z\u{a}linescu, C. }\emph{Convex analysis in
general vector spaces}. World Scientific Publishing Co., Inc., River Edge, NJ, 2002.

\bibitem {Zheng}\textsc{Zheng, X. Y. and Ng, K. F.} \emph{Subsmooth
semi-infinite and infinite optimization problems}. Math. Program. 134 (2012),
Ser. A, 365--393. 
\end{thebibliography}
\end{document}